\documentclass{article}      

\usepackage{amssymb}      
\usepackage{amsmath}      
\usepackage{amsthm}      
\usepackage{amsbsy}      
\usepackage{amscd}      
\usepackage[dvips]{graphics}
 \usepackage{mathrsfs}

\theoremstyle{plain}      
\newtheorem{theorem}{Theorem}[section]

\newtheorem{lemma}[theorem]{Lemma}      
\newtheorem{corollary}[theorem]{Corollary}      
\newtheorem{proposition}[theorem]{Proposition}

\newtheorem{definition}{Definition}[section]          
\theoremstyle{remark}      
\newtheorem{remark}[theorem]{Remark}  
\newtheorem{notation} [theorem]{Notation}

\newcommand{\Q}{\mathbb{Q}}      
      
\newcommand{\Z}{{\mathbb{Z}}}   
 \newcommand{\N}{{\mathbb{N}}}  
\newcommand{\R}{{\mathbb{R}}}

 \newcommand{\C}{{\cal C}^+(\T)}
\newcommand{\CR}{{\cal C}(\T)}
       
\newcommand{\B}{{{\cal B}}}

\newcommand{\T}{T^*}

\newcommand{\aps}{{{\alpha^{\star}}}} 
\newcommand{\als}{{\alpha^*}}      
\newcommand{\bps}{{{\beta^{\star}}}} 
\newcommand{\bs}{{\beta^*}}

\begin{document}

\date{\today}

\title{
The braided Ptolemy-Thompson group is finitely presented\footnote{     
This version: {\tt June 2005}.      
This preprint is available electronically at       
          \tt  http://www-fourier.ujf-grenoble.fr/\~{ }funar }}      
 \author{      
\begin{tabular}{cc}      
 Louis Funar &  Christophe Kapoudjian\\      
\small \em Institut Fourier BP 74, UMR 5582       
&\small \em Laboratoire Emile Picard, UMR 5580\\      
\small \em University of Grenoble I &\small \em University of Toulouse      
III\\      
\small \em 38402 Saint-Martin-d'H\`eres cedex, France      
&\small \em 31062 Toulouse cedex 4, France\\      
\small \em e-mail: {\tt funar@fourier.ujf-grenoble.fr}      
& \small \em e-mail: {\tt ckapoudj@picard.ups-tlse.fr} \\      
\end{tabular}      
}      

\maketitle

\begin{abstract} 
Pursueing our investigations on the relations between Thompson groups
and mapping class groups,  we introduce the group $T^*$ (and its 
further generalizations) which is an extension  of the Ptolemy-Thompson group $T$ by means of the 
full braid group $B_{\infty}$ on infinitely many strands. We prove that it is a finitely presented group with solvable
word problem, and give an explicit presentation of it.
\vspace{0.2cm} 
 
\noindent 2000 MSC Classification: 20 F 36, 20 F 38, 20 F 05, 57 M 07, 57 N 05  
 
\noindent Keywords: braid groups, mapping class groups, infinite surface, Thompson group.  
 
\end{abstract}

\section*{Introduction}      

In \cite{FK} we introduced a group ${\mathcal B}$ which is a universal 
mapping class group of genus zero. This is a mapping class 
group (in a restricted sense) of the Cantor surface of genus zero 
consisting of a 2-dimensional sphere with a Cantor set deleted. 
The main result of \cite{FK} 
is a finite presentation of  ${\mathcal B}$ coming from a cocompact 
action of the group on an explicit  complex, with small stabilizers. 
On the other hand  ${\mathcal B}$ may be algebraically described as an extension of the 
Thompson group $V$ by a pure mapping class group of the Cantor surface. 

\vspace{0.2cm} 
\noindent 
The aim of this paper is to pursue further  these lines of thought. 
One of the motivations of \cite{FK} was to prepare the setting for the finite 
presentability of a universal mapping class group of infinite genus, 
constructed in the same way as ${\mathcal B}$, but using a surface 
of infinite genus. The next step in completing this project is the 
knowledge of an intermediary group ${\mathcal B}^*$ that is associated  
to the infinite surface with handles surgered off and replaced by punctures.  
Our result in this direction is that this group is  also finitely presented. 

\vspace{0.2cm} 
\noindent 
However, there exists a simpler construction that presents the same 
features as  the shift from  ${\mathcal B}$ to ${\mathcal B}^*$.  
This construction starts from the Thompson group $T$ and yields an 
extension $T^*$ 
of $T$ by the group of braids $B_{\infty}$ on infinitely many strands. 
Notice that an extension $BV$ of the larger Thompson group $V$ by a pure braid group 
 on infinitely many strands has been recently considered 
by M. Brin and P. Dehornoy (\cite{Bri0,Bri,De1,De2}). It constitutes
the planar counterpart of the  group ${\mathcal B}$, in which the 
Cantor surface is replaced by a disk minus a Cantor set. As a matter of fact, we have observed that $BV$ is a subgroup of ${\mathcal B}$ (see \cite{FK}). 
Since
$BV$ is called the {\it braided Thompson group}, we hope to avoid any confusion by calling $T^{*}$ the 
{\it braided Ptolemy-Thompson group}, insisting on its relation with the Ptolemy-Thompson
group $T$ (see \cite{pe}). The group $T^{*}$ is essentially
 different from $BV$ (and ${\cal B}$), being an extension by the whole group of braids, and not only the pure braids.

\vspace{0.2cm} 
\noindent 
There is a group in the literature, indeed, to which $T^{*}$ should be compared. This is the somewhat mysterious acyclic
 extension 
considered earlier by Greenberg and Sergiescu (\cite{GS}). A recent construction of the Greenberg-Sergiescu extension has been given in \cite{KS}, in terms
of a mapping class group of an infinite surface with punctures. Actually, we define $T^*$ as a mapping class group
 of a planar surface $D^*$ obtained by thickening in the plane the binary tree. The elements of $T^*$ are 
 required to permute the punctures of $D^*$ (corresponding to the vertices of the tree) and to preserve at infinity a 
 certain rigid structure of the surface. In this respect, our construction is close to that of \cite{GS}. It is even simpler, but (and since) our
group $T^{*}$, contrary to the Greenberg-Sergiescu group, does not encode the discrete Godbillon-Vey class. The algebraic properties 
of the Greenberg-Sergiescu group are very poorly understood today. Our second motivation in studying $T^*$ is 
to get some insight which would enable us to understand the 
Greenberg-Sergiescu group. 

\vspace{0.2cm} 
\noindent 
Lastly, the construction of $T^*$, as  a mapping class group of a punctured infinite surface 
of infinite type, can be viewed also as an extension 
of the Dynnikov three pages representations (\cite{Dy}), where, for the 
first time, the infinite braid group $B_{\infty}$ was realized as the  
commutator of a finitely presented group.  Notice that this way 
of encoding links as three pages braids led Dynnikov to a purely algebraic 
algorithm for recognizing the unknot.

\vspace{0.2cm} 
\noindent 
Our main result is the following:
\begin{theorem}
The braided Ptolemy-Thompson group $T^*$ is a finitely presented group. 
\end{theorem}

\vspace{0.2cm} 
\noindent The bulk of the paper is devoted to the proof of this theorem, with an explicit presentation for $T^*$. We follow K. Brown's method (see \cite{br}), based on the Bass-Serre theory. It
 consists in building
a simply connected 2-dimensional complex on which $T^*$ acts cocompactly with finitely presented stabilizers. 
The complex is a kind of fibration over a (conveniently reduced) Hatcher-Thurston complex of the infinite surface. The latter is
a quotient of the Cayley complex of the Ptolemy-Thompson group $T$. The terminology used here for $T$ is expected
to stress on its relation with the Penner-Ptolemy groupoid (see \cite{pe}), from which its presentation actually derives.\\
 A similar construction was used in \cite{FK} to build up a complex for 
${\mathcal B}$,  but there, the problem was rather complicated 
because of the complexity of the Brown-Stein complex for $V$.\\ 
However, our $T^*$-complex has a new ingredient: the fibre of the fibration over the Hatcher-Thurston complex
is the Cayley complex
of the braid group $B_{\infty}$. A key tool is the use of a presentation of $B_{\infty}$ which exploits the
homogeneity of the tree associated to the infinite surface. It is precisely provided by a general theorem due to V. Sergiescu (\cite{Se}).
The remaining difficulty is that countably many relations of commutations between the braid generators occur in this presentation. 
Fortunately, the dependence between the braidings and the Thompson generators is so strong that we
manage to make use of the finiteness properties of $T$ in order to eliminate almost all the cumbersome commutation relations in $B_{\infty}$.

\vspace{0.2cm} 
\noindent  The same kind of methods also works for the group 
${\mathcal B}^*$ which is an extension of ${\mathcal B}$ by the infinite 
braid group of the Cantor surface. 
We state the result at the end of the paper.  

\vspace{0.2cm} 
\noindent We also introduce a group $T^{\star}$, a twin brother of $T^*$, which satisfies
the same finiteness property as $T^*$, but has only two generators, while $T^*$ has three. This provides the following new formulation
for the main theorem:
\begin{theorem}
The stable braid group $B_{\infty}$ embeds as a normal subgroup into the finitely presented group $T^{\star}$ with two generators.
\end{theorem}

\noindent The groups $T^{\star}$ and $T^*$ generalize the diagram picture groups considered by V. Guba and M. Sapir (\cite{GuS}) insofar as these 
are extensions by infinite permutation groups, rather than braid groups. Diagram groups are known to have very good
properties: they are $FP_{{\infty}}$, have solvable conjugacy and word problems.  We will prove that $T^{\star}$ and $T^*$ 
also satisfy the following:
\begin{theorem}
The groups $T^*$ and $T^{\star}$ have solvable word problems.
\end{theorem}

\vspace{0.2cm} 
\noindent {\bf Acknowledgements.} The authors are thankful to D. Calegari,  P. Dehornoy, M. Rubin, V. Sergiescu and 
B. Wiest, for comments and useful discussions.

\section{Infinite planar surfaces  and asymptotic mapping class groups}      

\subsection{Enhanced surfaces of infinite type} 

The surfaces below will be oriented and all      
homeomorphisms considered in the sequel will be      
orientation-preserving, unless the opposite is explicitly stated.   
\begin{definition}\label{ss}      
The ribbon tree $D$ is the planar surface obtained by   
thickening in the plane the infinite binary        
tree. We denote by $D^{*}$ (respectively,  $D^{\star}$) the ribbon tree  
with infinitely many punctures, one puncture for each vertex (respectively, each edge) of the 
tree.  
\end{definition}      
      
 \begin{figure} 
  \begin{center}   
\includegraphics{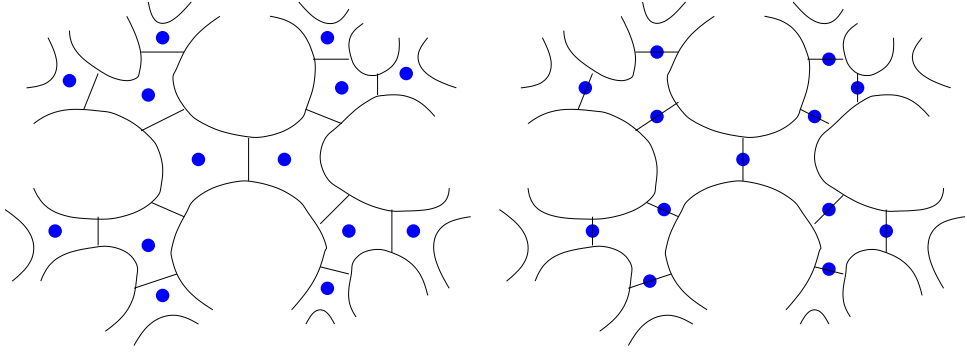}      
\caption{$D^*$ and $D^{\star}$ and with their canonical rigid structures}\label{tt}    
\end{center}  
\end{figure}

\begin{definition}
A {\em rigid structure} on $D$, $D^*$ or $D^{\star}$ is a decomposition  into
hexagons by means of a family of arcs whose 
endpoints are on the boundary of $D$. Each puncture lies inside a hexagon in the case of $D^*$, while each
 arc passes through a puncture in the case of $D^{\star}$. It is assumed that these arcs  
are pairwise non-homotopic in $D$, by homotopies keeping the boundary points of the arcs 
on the boundary of $D$. The choice of a rigid structure of reference is called the {\em canonical rigid structure}. 
The canonical rigid structure 
of the ribbon tree $D$ is such that each arc of this rigid structure crosses once and transversely a unique edge of the tree.
 The canonical 
rigid structures on $D^*$ and $D^{\star}$ are assumed to coincide with the canonical rigid structure of $D$ when forgetting the punctures. 
See Figure \ref{tt}. 
\end{definition}
 
\begin{notation}\label{not}  
 Let $^{\diamondsuit}$ stand for $^*$,  $^{\star}$  or the vacuum. The set of isotopy classes of
 rigid structures on $D^{\diamondsuit}$ will be denoted $\mathfrak{R}^{\diamondsuit}$. The canonical rigid structure of
 $D^{\diamondsuit}$ will be denoted $\mathfrak{r}^{\diamondsuit}$.
 \end{notation}
 
\subsection{Asymptotic mapping class groups}      

\begin{definition} 
\noindent 1. Let $D^{\diamondsuit}$ denote $D$,  $D^*$ or $D^{\star}$.  
A planar subsurface of $D^{\diamondsuit}$ is {\em admissible}  if  
it is a connected finite union of hexagons coming from the canonical rigid structure $\mathfrak{r}^{\diamondsuit}$.  
The {\em frontier} of an admissible surface is the  
union of the arcs contained in the boundary.\\
\noindent 2. Let $\varphi$ be a homeomorphism of $D^{\diamondsuit}$. One says that       
$\varphi$ is {\em asymptotically rigid} if the following conditions are      
fulfilled:      
\begin{itemize}      
\item There exists an admissible subsurface $\Sigma\subset 
  D^{\diamondsuit}$ such  that $\varphi(\Sigma)$ is also admissible.       
\item The complement $D^{\diamondsuit} -\Sigma$ is a union of $n$ infinite      
  surfaces. Then the restriction  
$\varphi: D^{\diamondsuit} -\Sigma\to D^{\diamondsuit}-\varphi(\Sigma)$      
 is {\em rigid}, meaning that it respects the rigid structures in the      
 complements of the compact subsurfaces, mapping hexagons into hexagons.       
  Such a surface $\Sigma$ is called a {\em support} for $\varphi$.      
\end{itemize}      
One denotes by $T$, $T^*$ and $T^{\star}$ the group of isotopy classes of asymptotically rigid  homeomorphisms of $D$,
 $D^*$ or $D^{\star}$, respectively.
    
\end{definition} 
\begin{remark} 
There exists a cyclic order on the frontier arcs of an  
admissible subsurface induced by the planarity. An asymptotically 
rigid homeomorphism necessarily preserves the cyclic order  
of the frontier for any admissible subsurface.  

\end{remark} 
 
\subsection{Ptolemy-Thompson  group $T$ as a mapping class group}
The mapping class group $T$ is isomorphic to the Thompson group which is commonly denoted $T$. This fact has been widely developed in \cite{KS} and \cite{FK}.
We consider the following elements of $T$, defined as mapping classes 
of asymptotically rigid homeomorphisms: 
\begin{itemize} 
\item The support of the element $\beta$ is the central hexagon.  
Further,  $\beta$ acts as the counterclockwise rotation of order  
three which permutes the three branches 
of the ribbon tree issued from the hexagon.  

\begin{center} 
\includegraphics{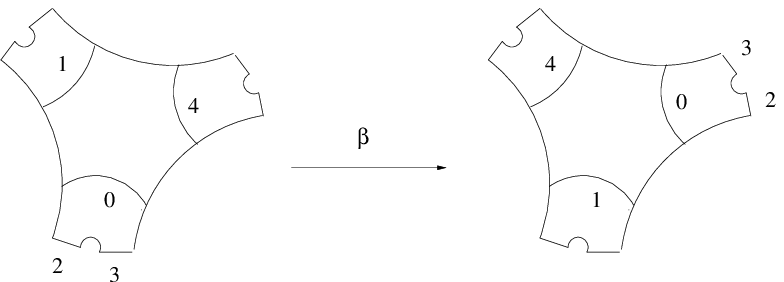}      
\end{center}
\item  The support of $\alpha$  is the union of two adjacent hexagons, 
one of them being the support of $\beta$.  Then $\alpha$  
rotates counterclockwise the support of angle $\frac{\pi}{2}$, by  
permuting the four branches of the ribbon tree
issued from the support.  

\begin{center} 
\includegraphics{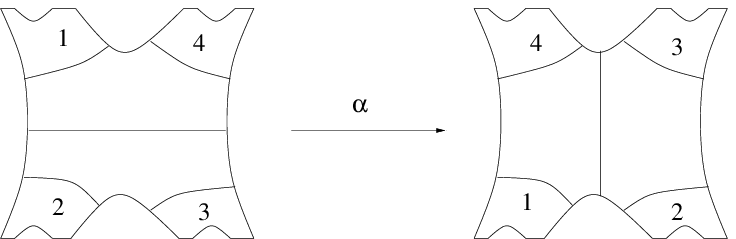}      
\end{center}
\end{itemize} 

\begin{proposition}\label{presT}
The group $T$ has the following presentation with generators 
$\alpha$ and $\beta$ and relations 
\[ \alpha^4=\beta^3=1\]
\[ [\beta\alpha\beta, \alpha^2\beta\alpha\beta\alpha^2]=1\]
\[  [\beta\alpha\beta, \alpha^2\beta^2\alpha^2\beta\alpha\beta\alpha^2
\beta\alpha^2]=1\]
\[ (\beta\alpha)^5=1\]
\end{proposition}
\begin{proof} This result is due to Lochak and Schneps 
(see \cite{LS}), but notice that there is a typo in their statement, 
which is corrected above. 
\end{proof}

\begin{remark}
If we set $A=\beta\alpha^2$, $B=\beta^2\alpha$ and 
$C=\beta^2$ then we obtain the generators $A,B,C$ of the group $T$, 
considered in \cite{CFP}. Then the  first two relations above 
are equivalent to 
\[ [AB^{-1}, A^{-1}B A]=1, \,\;\; [AB^{-1}, A^{-2}B A^2]=1\]
The presentation of $T$ in terms of the generators 
$A,B,C$ consists  of the two relations above with  four more 
relations to be added: 
\[ C^3=1, \;\, C=B A^{-1}CB,\,\; CA=(A^{-1}CB)^2,\,\; (A^{-1}CB)(A^{-1}BA)=B
(A^{-2}CB^2)\]
\end{remark}
\begin{remark}
The subgroup of $T$ generated by the elements $A$ and $B$ 
is the Thompson group $F$, as it is obvious from \cite{CFP}. 
Moreover, the group $F$ has the presentation  
\[ F=\langle A, B; [AB^{-1}, A^{-1}B A]=1, \,\,  [AB^{-1}, A^{-2}B A^2]=1\rangle \]
Consequently, the inclusion $F\to T$  sends $A$ to $A$ and $B$ to $B$. 
\end{remark}

\section{The braided Ptolemy-Thompson groups $T^*$ and $T^{\star}$}
\subsection{$T^*$ and $T^{\star}$ as extensions of the Thompson group $T$} 
\noindent We write $D$ (respectively $D^*$ and $D^{\star}$) as  
an ascending union $\cup_n D_n$, where $D_0$ is the support of 
$\beta$, and $D_{n+1}$ is obtained by adding a new hexagon to $D_n$   
along each component of the frontier.\\   
Let the symbol $^{\diamondsuit}$ denote either $^*$ or $^{\star}$. The Artin braid groups $B[D^{\diamondsuit}_{n}]$ associated to the punctures  
on $D^{\diamondsuit}_{n}$ form an inductive system induced by the inclusions  
 $D_n^{\diamondsuit}\subset D_{n+1}^{\diamondsuit}$, whose limit  
$B[D^{\diamondsuit}]=\lim_{n\to\infty}B[D_n^{\diamondsuit}]$ can be identified  
with the group of compactly supported braids on $D$, where the base points  
of the strands are the punctures of $D^{\diamondsuit}$. 

\begin{remark}
Let the symbol $^{\diamondsuit}$ denote either $^*$ or $^{\star}$. The group $B[D^{\diamondsuit}]$ is isomorphic to the
stable braid group
$B_{\infty}$, where $B_{\infty}$ is the inductive
limit coming from the inclusions $\sigma_{i}\in B_{n}\mapsto \sigma_{i}\in B_{n+1}$, where $\sigma_{i}$ ($1\leq i\leq n-1$)
denotes a standard Artin generator. This can be proven by observing first that the embedding of $D^{\diamondsuit}$ into 
the Euclidean plane $P$ induces 
an isomorphism of $B[D^{\diamondsuit}]$ with the group of isotopy classes of compactly supported homeomorphisms
 of the punctured plane. Since the set of punctures is discrete in $P$, one may construct a homeomorphism of $P$ which
 maps the punctures on the points of coordinates $(i,0)$, $i\in\N^*$ (after the choice of a framing).  By conjugation,  this
  homeomorphism
 induces an isomorphism between $B[D^{\diamondsuit}]$ and $B_{\infty}$.
\end{remark}

\begin{proposition}\label{plan} 
Let the symbol $^{\diamondsuit}$ denote either $^*$ or $^{\star}$. We have an exact sequence:  
\[  
1 \to B[D^{\diamondsuit}] \to T^{\diamondsuit}  \to  T          \to   1.  
\] 
\end{proposition} 
\begin{proof} 
Thompson's group $T$ is viewed here as the group of asymptotically rigid homeomorphisms of $D$ 
(without punctures) 
up to isotopy. Thus, the epimorphism $ T^{\diamondsuit}  \to  T$ is induced by forgetting the punctures. 
Now let $\varphi$ be an asymptotically rigid  homeomorphism of 
$D^{\diamondsuit}$ whose image in $T$ is trivial. This implies that outside an admissible subsurface, $\varphi$ is
isotopic to identity. Without changing the class of $\varphi$ in $T^{\diamondsuit}$, we may assume that outside this subsurface,
 $\varphi$ {\it is}  identity. Therefore, there exists a compactly supported  
isotopy $\varphi_t$ among homeomorphisms of $D$ which  
joins $\varphi$ to identity, whose support is an admissible subsurface. Further the class of the homeomorphism $\varphi$ on the punctured support is  
completely encoded by a braid, and a picture showing the trajectory  
of the punctures during the isotopy. Then the class of $\varphi$  
is the image of a braid from a some suitable $B[D_n^{\diamondsuit}]$
into $T^{\diamondsuit}$. This means that the kernel above is $B[D^{\diamondsuit}]$. 
\end{proof}
 
\subsection{$T^{\star}$ is generated by two elements}\label{2el}
Thompson's groups and their generalizations considered by 
Higman are generated by two elements
(\cite{Mas}). It is known that mapping class groups
of closed surfaces of genus at least one are also 
generated by two elements. We will prove here that 
the same holds for the group $T^{\star}$.\\
  
\noindent{\bf Specific elements}\\
\noindent Let us consider the following elements of $T^*$ and $T^{\star}$: 
\begin{itemize} 
\item The support of the element $\beta^{\star}$ of $T^{\star}$ (resp. $\beta^*$ of $T^*$) is the central hexagon.  
Further $\beta^{\star}$ and $\beta^*$ act as the counterclockwise rotation of order  
three which permutes cyclically  
the punctures.  One has ${\beta^{\star}}^{3}=1$ and ${\beta^*}^{3}=1$.

\begin{center} 
\begin{picture}(0,0)%
\includegraphics{beta.pstex}%
\end{picture}%
\setlength{\unitlength}{987sp}%
\begingroup\makeatletter\ifx\SetFigFont\undefined%
\gdef\SetFigFont#1#2#3#4#5{%
  \reset@font\fontsize{#1}{#2pt}%
  \fontfamily{#3}\fontseries{#4}\fontshape{#5}%
  \selectfont}%
\fi\endgroup%
\begin{picture}(27174,4512)(2764,-4411)
\put(23626,-1711){\makebox(0,0)[lb]{\smash{{\SetFigFont{6}{7.2}{\familydefault}{\mddefault}{\updefault}{$\beta^*$}%
}}}}
\put(8626,-1786){\makebox(0,0)[lb]{\smash{{\SetFigFont{6}{7.2}{\familydefault}{\mddefault}{\updefault}{$\beta^{\star}$}%
}}}}
\end{picture}%
\
\end{center}
\item  The support of the element $\alpha^{\star}$ of $T^{\star}$ (resp. $\alpha^*$ of $T^*$) is the union of two adjacent hexagons, 
one of them being the support of $\beta^{\star}$ and $\beta^*$.  
Then $\alpha^{\star}$  (resp. $\alpha^*$)
rotates counterclockwise the support 
of angle $\frac{\pi}{2}$, by  
keeping fixed the central puncture (resp. the two punctures of the adjacent hexagons). One has  ${\alpha^{\star}}^4=1$
 while ${\alpha^*}^4=\sigma^2$, where $\sigma$ denotes the braid that permutes the puncture 0 and 3, see Definition
  \ref{defbraid} below.

\begin{center} 
\begin{picture}(0,0)%
\includegraphics{alpha.pstex}%
\end{picture}%
\setlength{\unitlength}{1184sp}%
\begingroup\makeatletter\ifx\SetFigFont\undefined%
\gdef\SetFigFont#1#2#3#4#5{%
  \reset@font\fontsize{#1}{#2pt}%
  \fontfamily{#3}\fontseries{#4}\fontshape{#5}%
  \selectfont}%
\fi\endgroup%
\begin{picture}(25149,3699)(2689,-10723)
\put(8326,-8686){\makebox(0,0)[lb]{\smash{{\SetFigFont{7}{8.4}{\familydefault}{\mddefault}{\updefault}{$\alpha^{\star}$}%
}}}}
\put(21826,-8686){\makebox(0,0)[lb]{\smash{{\SetFigFont{7}{8.4}{\familydefault}{\mddefault}{\updefault}{$\alpha^{\star}$}%
}}}}
\end{picture}%
\end{center}
\end{itemize} 

\begin{definition}\label{defbraid}
Let $e$ be a simple arc in $D^*$ or $D^{\star}$ which connects two punctures. We associate 
a braiding $\sigma_e\in B_{\infty}$ to $e$ by considering the homeomorphism
that moves clockwise the punctures at the 
endpoints of the edge $e$ in a small neighborhood of the edge,  
in order to interchange their positions. This means that, if $\gamma$ is an arc transverse to $e$, then the braiding $\sigma_{e}$
moves $\gamma$ on the left when it approaches $e$. Such a braiding will be called {\it positive}, while $\sigma_{e}^{-1}$ is
{\it negative}.
\end{definition}

\begin{center} 
\includegraphics{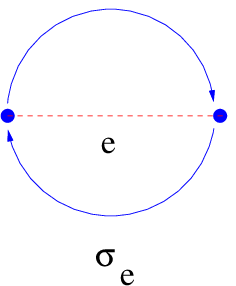}      
\end{center}

\begin{remark}\label{PSL}
The subgroup of $T^{\star}$ generated by ${\alpha^{\star}}^2$ and $\beta^{\star}$ is isomorphic to $PSL(2,\Z)$, viewed
as the group of orientation-preserving automorphisms of the binary planar tree of the ribbon surface $D$. In the same way, the
 subgroup of $T^*$ generated by
 $a={\alpha^*}^2\sigma^{-1}$ (which is of order 2) and 
$\beta^*$ is isomorphic to $PSL(2,\Z)$.
\end{remark}

\begin{theorem} 
$T^{\star}$ is generated by $\aps$ and $\bps$. 
\end{theorem} 
\begin{proof}

\noindent 
Let us denote by $T'$ the subgroup  of $T^{\star}$ generated by $\aps,\bps$.  
\begin{lemma}
The restriction to $T'$ of the projection map $T^{\star}\to T$ is surjective. 
\end{lemma}
\begin{proof}
This maps sends $\aps$ to $\alpha$ and $\bps$ to $\beta$. According to Proposition \ref{presT},  $\alpha$ and $\beta$ generate 
$T$, thus the claim. 
\end{proof}

\vspace{0.2cm} 
\noindent 
Let now $\iota:B_{\infty}\to T^{\star}$ be the natural inclusion. 
Since $T^{\star}$ is an extension of $T$ by $B_{\infty}$, it suffices now to show 
that $\iota(B_{\infty})\subset T'$. This will be done in two steps. 
First we show that a specific braid  generator lies in $T'$ and next we use 
the conjugation action to prove that all braid generators lay within $T'$. 
The first steps proceeds as follows. 

\begin{lemma}\label{sigmapenta}
The braid generator $\sigma_{[04]}$ associated to the edge joining the 
punctures numbered $0$ and $4$ has the image 
\[ \iota(\sigma_{[04]})=(\bps\aps)^5\]
\end{lemma} 
\begin{proof}
An explicit picture calculations shows that the action of 
$(\bps\aps)^5$ on the  standard rigid structure  
of the ribbon tree is the following one: 

\begin{center}     
\begin{picture}(0,0)%
\includegraphics{sigma.pstex}%
\end{picture}%
\setlength{\unitlength}{1381sp}%
\begingroup\makeatletter\ifx\SetFigFont\undefined%
\gdef\SetFigFont#1#2#3#4#5{%
  \reset@font\fontsize{#1}{#2pt}%
  \fontfamily{#3}\fontseries{#4}\fontshape{#5}%
  \selectfont}%
\fi\endgroup%
\begin{picture}(12174,3774)(2689,-10723)
\put(8026,-8536){\makebox(0,0)[lb]{\smash{{\SetFigFont{8}{9.6}{\familydefault}{\mddefault}{\updefault}{$(\beta^{\star}\alpha^{\star})^5$}%
}}}}
\end{picture}%

\end{center}

\noindent In particular, this can be identified with the action 
of the mapping class $\sigma_{[04]}$. 
Then the action of $(\bps\aps)^5$ coincides with the natural action  
of $\sigma_{[04]}\in B_{\infty}$ 
on the arcs in the punctured surface. In meantime the  
configuration of arcs  coming from a rigid structure 
(up to isotopy) uniquely determines   
the element of $B_{\infty}$, and so $\iota((\bps\aps)^5)=
\sigma_{[04]}$. 
\end{proof}

The end of the proof is now as follows. For each hexagon of $D^{\star}$, consider the three 
arcs which connect the punctures to each other, and  intersect only at the punctures. Let $\cal E$ be the set of all such arcs associated
to $D^{\star}$. The subgroup of $T^{\star}$ generated by $\aps^2$ and $\bps$ acts transitively on $\cal E$ as the
 group $PSL(2,\Z)$. Therefore, for each $e\in\cal E$, there exists a word $w$ in $\aps^2$ and $\bps$ such that
  $w(e_{[04]})=e$. Then $w\sigma_{[04]}w^{-1}=\sigma_{e}$. Consequently, each $\sigma_{e}$ belongs to $T'$. Since
  the group $B_{\infty}$ is generated by the braidings $\sigma_{e}$ when $e$ runs over $\cal E$, $T'$
   contains $\iota(B_{\infty})$.
\end{proof}

\begin{figure} 
\begin{center}     
\includegraphics{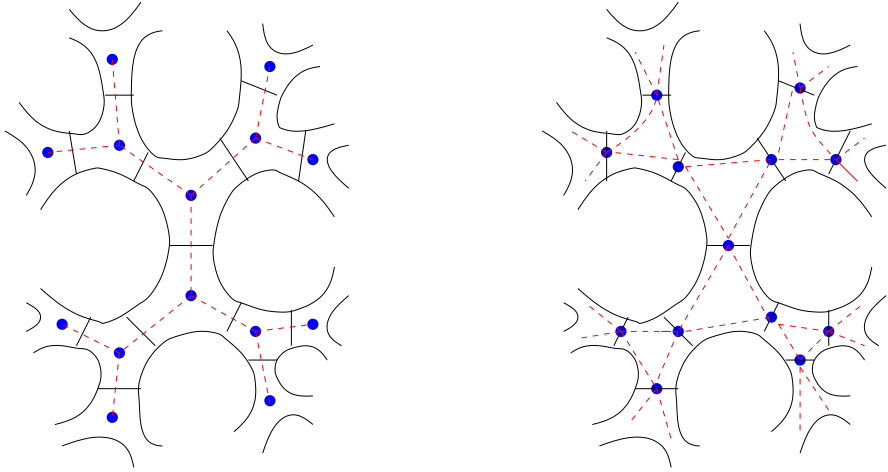}      
\end{center}
\caption{Tree of $D^*$ and Graph of $D^{\star}$}\label{arbres}      
\end{figure}

\begin{remark}\label{graph}
The union of all edges of $\cal E$ is a graph, which is dual to the binary tree (of $D$ or $D^*$). It will be called {\it the graph of 
$D^{\star}$}, see Figure \ref{arbres}.
A general theorem due to V. Sergiescu 
(\cite{Se})
implies that $B_{\infty}$ is generated by $\{\sigma_{e},\;e\in{\cal E}\}$. The relations holding between these
 generators are explicited in \cite{Se}. This approach was later generalized by Birman, Ko and Lee (\cite{BKL}). 
\end{remark}

\subsection{Solution of the word problem}
\begin{proposition}\label{soluble}
The groups $T^*$ and $T^{\star}$ have solvable word problem. 
\end{proposition}
\begin{proof}
The proofs for $T^*$ and $T^{\star}$ being the same, we only consider the case of $T^{\star}$.
Consider a word  $w$ in the free group generated by the 
letters $\aps,\bps$, of length $|w|=n$. Our aim is to 
find an algorithm to decide whether the element 
represented by $w$ in $T^{\star}$ is trivial or not. When the former 
alternative holds true, we will write this as $w =_{_{T^{\star}}}\,1$. 
We will also denote by $[w]$ the class of the word $w$ in $T^{\star}$. 

\vspace{0.2cm}
\noindent Recall that we denoted by $D_0^*$ the support 
of $\bps$ (the central 
hexagon) and that $D_{n+1}^{\star}$ is the result of 
adding a  new hexagon along 
each boundary component of $D_n^{\star}$. Thus $D^{\star}_n$ has 
$3\cdot 2^{n+1} -3$ punctures, 
$3\cdot 2^n-2$ hexagons and $3\cdot 2^{n}$ boundary arcs. 
We will say that the boundary hexagons of 
$D^{\star}_n$  - i.e. those  which don't belong to $D^{\star}_{n-1}$ - 
are at distance $n$ from $D^{\star}_0$. The distance function 
between hexagons is induced by the  usual distance at the 
level of the dual tree. 

\vspace{0.2cm}
\noindent The main ingredient is the following  {\em a priori} 
characterization of the size  of supports:

\begin{lemma}\label{supp} 
If $|w|=n$ then the support of $[w]$ is contained in $D_n^{\star}$. 
This means that there exists an admissible subsurface 
$\Sigma\subset D_n^{\star}$ such that $[w](\Sigma)$ is also admissible 
and $[w]:D^{\star}-\Sigma\to D^{\star}-[w](\Sigma)$ is rigid. 
Moreover, $[w](\Sigma)\subset D^{\star}_{n}$. 
\end{lemma}
\begin{proof}
First, we have $\bps^m(D^{\star}_n)=D^{\star}_n$ and  
$\aps^{m}(D^{\star}_n)\subset D^{\star}_{n+1}$ for any $m$. Recall that 
$\aps$ and $\bps$ are of finite order. 

\vspace{0.2cm}
\noindent
The claim holds true trivially for $n=1$. We use induction on $n$. 
Thus, if $|w|=n$ and $[w]$ has support $\Sigma$ as claimed, then 
$[\bps^m w]$ has support 
$\Sigma\cup \bps^{-m}(\Sigma)\subset D^{\star}_n$. 
Moreover, $[\aps^m w]$ has support 
$\Sigma\cup \aps^{-m}(\Sigma)\subset  D^{\star}_{n+1}$. 
The same argument works for $[w](\Sigma)$. 
This completes the induction step. 
\end{proof}
\begin{remark}
Actually, we proved that the support of $[w]$ is contained in  
$D^{\star}_k$, where $k$ is bounded from above by 
the number of distinct factors $\aps^m$, $\bps^m$ in 
the word $w$.  
\end{remark}

\noindent 
Consider now the associated word $\overline{w}$ obtained from 
$w$ by replacing $\aps,\bps$ by $\alpha,\beta$. Then  
$\overline{w}$ can be seen as a word representing 
the element $\pi([w])\in T$. 

\vspace{0.2cm}
\noindent 
The group $T$ is finitely presented and simple. 
Therefore by a well-known result $T$ has 
soluble word problem. Thus there exists an algorithm which decides 
whether $\overline{w}=_{_{T}}\,1$. 

\vspace{0.2cm}
\noindent
Assume that $\overline{w}=_{_{T}}\,1$. This implies that 
the element $[w]\in B_{\infty}\subset T^*$. According to the 
previous lemma $[w]\in B(D^{\star}_n)\subset B_{\infty}$, where 
$B(D^{\star}_n)$ is the braid group of the $D^{\star}_n$ and thus it suffices to 
decide whether the image of $w$ is trivial in $B(D^{\star}_n)$. 
However, $w$ is {\em not} given as a word in the  
generators of the braid group, but as a word in $\aps,\bps$ 
which - one knows  that - it {\em can} be rewritten as a word 
in the generators. Thus we have to overcome the difficulties  
concerning the rewriting process.

\vspace{0.2cm}
\noindent
Actually, there exists an algorithm which rewrites  
$[w]$ as  a word $\sigma_{i_1}^{m_1}\sigma_{i_2}^{m_2}\cdots 
\sigma_{i_k}^{m_k}$, where $\sigma_i$ are the standard 
generators of $B[D^{\star}_n]$, coming from twists of endpoints of the 
edges of a maximal tree. In particular,  
$[w]=_{_{T^{\star}}}\,\sigma_{i_1}^{m_1}\sigma_{i_2}^{m_2}\cdots \sigma_{i_k}^{m_k}$. Let us assume that. The
 braid group is automatic and hence 
it has solvable word problem. This implies that  
it can be explicitly checked 
whether the word $\sigma_{i_1}^{m_1}\sigma_{i_2}^{m_2}\cdots 
\sigma_{i_k}^{m_k}$ represents the identity in $B[D^{\star}_n]$. 
This is equivalent to $w=_{_{T^{\star}}}\, 1$. 
Notice that the result of this test is independent on the 
particular word  $\sigma_{i_1}^{m_1}\sigma_{i_2}^{m_2}\cdots 
\sigma_{i_k}^{m_k}$ chosen above.

\vspace{0.2cm}
\noindent
Here is an algorithm which permits to overcome the rewriting
procedure. The previous lemma shows that it suffices to understand 
the action of $[w]$ on $D^{\star}_n$, since $[w]$ is rigid on the 
complementary of $D^{\star}_n$ and thus identity. This implies that 
$[w]$ can be viewed as a mapping class (of a homeomorphism) of 
$D^{\star}_n$.  The  punctured disk $D^{\star}_n$ 
is decomposed into hexagons by means of some arcs passing thru 
the punctures, which form the rigid structure. The mapping class 
$[w]\in B[D^{\star}_n]$ is completely determined by the image of 
the arcs (or the rigid structure) up to isotopy. 
Furthermore, $w=_{_{T^{\star}}}\;1$ if and only if the image rigid 
structure is isotopic to the  initial rigid structure. 
Since the arcs are disjoint this is equivalent to saying that 
the image by $[w]$ of each arc is isotopic to itself. 
Notice that the isotopy can be supposed to be 
supported on $D^{\star}_n$, since $[w]$ is the identity outside $D^{\star}_n$.  

\vspace{0.2cm}
\noindent
The  first step is to understand what is needed for obtaining 
the image of the arc $\gamma\subset D^{\star}_n$ under the action 
of $[w]$. Let the input be the word 
$w=w_1w_2\cdots w_n$, where $w_j$ are  among the 
letters $\aps,\bps$ or their inverses. We know how to draw 
the image of an arc in $D^{\star}_n$ under each of the transformations 
$\aps, \bps$. However it may happen that the image of an 
arc $\gamma$ from $D^{\star}_n$ under $\aps$ be outside $D^{\star}_n$. 
Set then  $\gamma_1=w_n(\gamma)$ and   
$\gamma_{j+1}=w_{n-j}(\gamma_j)$, for $n-1\geq j\geq 1$. 
An application of the lemma \ref{supp} shows that 
if $\gamma$ is contained in $D^{\star}_n$, then $\gamma_j\in D^{\star}_{2n}$, 
for all $j\leq n$. Moreover, we have  $[w](\gamma)=\gamma_n$.  
In order to find the images of the arcs from $D^{\star}_n$ it 
suffices to understand the restriction of the action to 
$D^{\star}_{2n}$. Define now the restricted action of $\aps$ and $\bps$ 
on arcs of $D^{\star}_{2n}$  as follows: 
\[ \aps[2n] (\gamma)= \left\{\begin{array}{ll}
\aps(\gamma), & {\rm if }\; \aps(\gamma)\subset D^{\star}_{2n}\\
\emptyset, & {\rm otherwise}
\end{array}
\right.\]
\[ \bps[2n] (\gamma)= \left\{\begin{array}{ll}
\bps(\gamma), & {\rm if }\; \bps(\gamma)\subset D^{\star}_{2n}\\
\emptyset, & {\rm otherwise}
\end{array}
\right.\]
We define recurrently $w[2n](\gamma)=
w_1[2n](w_2[2n](\cdots (w_n[2n](\gamma))...)$. 
The previous arguments show
that $w[2n](\gamma)=[w](\gamma)$, when $\gamma\subset D^{\star}_n$.  
Thus, by restricting the original action 
on $D^{\star}_{2n}$ we get finiteness on one hand, and on the other 
hand  we are still able to recover the action of 
$[w]$ on $D^{\star}_n$ and so on all  of $D^{\star}$.

\vspace{0.2cm}
\noindent  
Recall now that simple curves on a surface are isotopic if and 
only if they are homotopic. Thus the problem on whether the 
the image by $[w]$ of the standard rigid structure
is isotopic to the former rigid structure is essentially 
an algebraic problem.
Fix a base point $*$ in $D^{\star}_n\subset D^{\star}$ and lifts of the 
mapping classes $\aps, \bps$ which preserve this base point. 
The transformations $\aps[2n], \bps[2n]$ act on loops based at $*$ 
and there are induced  explicit endomorphisms of free groups 
\[\aps[2n]: \pi_1(D^{\star}_{2n})\to \pi_1(D^{\star}_{2n})\] 
by setting $\aps[2n]([\gamma])=0$ if the image of $\gamma$ 
is not contained within $D^{\star}_{2n}$ (and similarly for $\bps[2n])$. 
The exact form of these endomorphisms can be read quickly by 
labeling the punctures, taking loops encircling the punctures 
as generators and identifying $\aps,\bps$ with 
planar rotations, but their exact form won't matter in the sequel. 
Define further  the restricted operators $w[2n](\gamma)=
w_1[2n](w_2[2n](\cdots (w_n[2n](\gamma))...)$.
These endomorphisms are {\em outer endomorphisms} which are 
well-defined only up to inner automorphisms.   
Consider now a basis $\gamma_1,\ldots,\gamma_N$ of 
$\pi_1(D^{\star}_n)$. Then the image of each element can be computed  
by using  only the restricted operators as follows: 
$[w](\gamma_j)=w[2n](\gamma_j)$. 
Remark that these formulas define an outer automorphism 
$[w]=w[2n]$ of $\pi_1(D^{\star}_n)$. 
Eventually, the test 
$[w]=_{{T^{\star}}}\,1$ is equivalent to checking whether 
this outer automorphism is trivial. Thus we have to check 
whether there exists $c\in \pi_1(D^{\star}_n)$ so that 
\[ (w[2n](\gamma_j))_{1\leq j\leq N} =_{_{\pi_1(D^{\star}_n)}} 
(c\gamma_jc^{-1})_{1\leq j\leq N} \]
This is equivalent to solving the  generalized 
conjugacy problem in the 
free group $\pi_1(D^{\star}_n)$. 
Now, the  generalized conjugacy 
problem is algorithmically solvable for such a group. 
This holds more generally for  biautomatic groups, as it follows 
immediately from the solution of the usual conjugacy problem (see  
\cite{ECHLPT}  Theorem 2.57 p.59, and \cite{GeSh}). 
This holds true also for groups satisfying the small 
cancellation conditions  C(6), $C(4)-T(4)$, 
$C(3)-T(6)$ (see \cite{Bez}). 
\end{proof}

\begin{remark}
The complexity of this algorithm is exponential 
(about ${\mathcal O}(n^d2^{n})$). In fact we can solve the problem 
$\overline{w}=_{_{T}}\, 1$ in about ${\mathcal O}(n^{14})$ steps
since the Dehn function of $T$ is bounded by $n^7$ 
(see \cite{Guba}). 
Next we can compute the  action of the 
restricted operators $w[2n](\gamma_k)$ in $n$ steps.
The length of the word $w[2n](\gamma_k)$ is ${\mathcal O}(n)$ and 
the conjugacy in the free group can be checked in a number 
of steps depending polynomially on the length. However the rank 
of the free group $\pi_1(D^{\star}_n)$ is $3\cdot 2^n$ and thus 
we have at least that amount of conjugacy tests.
\end{remark}

\begin{remark}
The solvability of the word problem is a strong indication that the 
group has nice properties from algorithmic viewpoint. Notice however 
that A. Yu. Olshanskii constructed  infinitely presented 
2-generator groups which have solvable word and 
conjugation problem. Thus the result above does not imply that 
$T^{\star}$ is finitely presented.  
\end{remark}

\subsection{The abelianization $H_1(T^{\star})$} 
\begin{proposition}
We have $H_1(T^{\star})=\Z/12\Z$. 
\end{proposition}
\begin{proof}
We know that $T^{\star}$ fits into an exact sequence:  
\[ 1\to B_{\infty}\to T^{\star}\to T\to 1\] 
The groups $B_{\infty}$ is the group of braids associated to the  
punctures of $D^{\star}$.  We will consider the Sergiescu presentation 
of $B_{\infty}$ associated to the graph of $D^{\star}$, see Remark \ref{graph}. Hall's lemma 
provides  an infinite presentation of $T^{\star}$:
the generators $\aps$ and $\bps$ satisfy  the following 
relations, as lifts of relations in $T$: 
\[ \aps^4=\bps^3=1\]
\[ [\bps\aps\bps, \aps^2\bps\aps\bps\aps^2]=1\]
\[  [\bps\aps\bps, \aps^2\bps^2\aps^2\bps\aps\bps\aps^2
\bps\aps^2]=1\]
\[ (\bps\aps)^5=\sigma_{[04]}\]
where $\sigma_{[04]}$ is the braid generator considered above. 
All relations involving the braids are coming from Sergiescu's 
relations above, while the remaining relations in $T^{\star}$ 
are conjugacy relations stating that $B_{\infty}$ is normal inside $T^{\star}$. 
Thus the 
abelianization $H_1(T^*)$ is generated by the classes of $\aps$ and $\bps$ 
which are constrained to be of order 4 and 3, respectively. 
The claim follows. 
\end{proof}

\subsection{Dynnikov's group as a mapping class group}
Dynnikov (\cite{Dy}) considered earlier a simpler finitely presented group 
which contains $B_{\infty}$, leading to the three page encoding of 
knots and links.  However, his example fits well in our more 
general framework. 
We will explicitly explain this on a specific example. 

\vspace{0.2cm}
\noindent Instead of considering  the ribbon tree obtained from the binary 
tree we consider the ribbon $Y$, where $Y$ denotes the wedge of three 
half-lines in the plane.  The group of automorphisms of $Y$ is then 
$\Z/3\Z$. Consider next $Y^*$ which is the ribbon $Y$  punctured 
at  the set of points of integer coordinates on each 
half-line, where the origin is assumed to be $0$. In the original definition
of Dynnikov the origin was not among the punctures, thus the group obtained below will be different from his group described in \cite{Dy}. 
There are three families of arcs, each one consisting 
of parallel arcs passing through the 
punctures of a half-line (excepting the central puncture) 
and connecting two boundary components. 
The surface $Y^*$ is then divided by these arcs 
into one hexagon  containing the  central puncture and infinitely many 
squares along the half-lines. 
One defines the admissible subsurfaces of $Y^*$  to be those 
hexagons determined by 
three arcs from three different families, and thus  
containing the central one. By analogy with the definitions of 1.1 and 1.2, we introduce:

\begin{definition} 
The group  $AS(Y^*)$ is the group of asymptotically rigid homeomorphisms of $Y^*$ 
up to isotopy.
\end{definition} 

We will consider next the 
subgroup $AS_{\partial}(Y^*)$ generated by  those homeomorphisms 
which are end preserving i.e. inducing a trivial automorphism of $Y$.
Alternatively, the homeomorphisms should send the arcs of some
support hexagons into corresponding arcs of the same family.

\begin{center}
\includegraphics{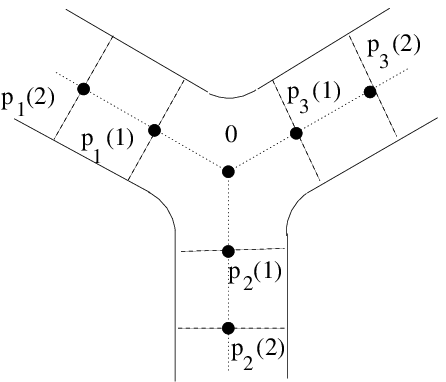}      
\end{center}

\begin{proposition}
There is an exact sequence 
\[ 1\to B_{\infty}\to AS_{\partial}(Y^*)\to \Z^2\to 1\]
where $B_{\infty}=\lim_{n\to\infty} B_{3n+1}$ is the limit of the braid 
groups of an exhausting sequence of subsurfaces of $Y^*$. 
\end{proposition}
\begin{proof} A mapping class $\varphi\in  AS_{\partial}(Y^*)$ 
sends a  support hexagon into another support hexagon, 
by translating the arc on the half-line 
$l_j$ of $n_j$ units towards the center. 
Since the support hexagons should contain 
the same number of punctures we have $n_1+n_2+n_3=0$.  
The map sending $\varphi$ to $(n_1,n_2,n_3)$ is a surjection onto 
$\Z^2$. The rest of the proof is immediate. 
\end{proof}

\vspace{0.2cm}
\noindent
Let the line $l_j$ be punctured along the points $p_j(i)$ at distance 
$i$ from the origin. Consider the  mapping class of the 
homeomorphism $d_j$ which translates all punctures of the line 
$l_{j-1} \cup l_{j+1}$ one unit in the counterclockwise direction, 
as in the figure below: 

\begin{center}
\includegraphics{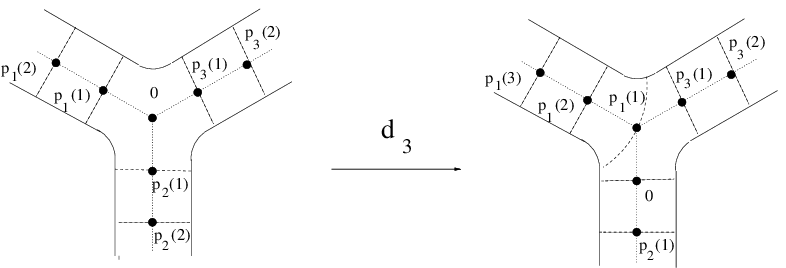}      
\end{center}

\begin{proposition}
The group  $AS_{\partial}(Y^*)$ is generated by the $d_1,d_2,d_3$ and admits 
the presentation given by the following relations: 
\[ d_1d_2d_3=1\]
\[ u_3u_2u_1u_3=u_1u_3u_2u_1=u_2u_1u_3u_2, 
{\rm where }~ u_i=d_{i+1}d_id_{i-1}\]
\[ u_1u_2u_1=u_2u_1u_2, u_3u_2u_3=u_2u_3u_2, u_1u_3u_1=u_3u_1u_3\]
\[ d_2^{-1}u_3d_2 u_2= u_2 d_2^{-1}u_3d_2, 
d_3^{-1}u_1d_3 u_3= u_3 d_3^{-1}u_1d_3, 
d_1^{-1}u_2d_1 u_1= u_1 d_1^{-1}u_2d_1\]
\[ d_2^{-1}u_3d_2 u_1= u_1 d_2^{-1}u_3d_2, 
d_3^{-1}u_1d_3 u_2= u_2 d_3^{-1}u_1d_3, 
d_1^{-1}u_2d_1 u_3= u_3 d_1^{-1}u_2d_1\]
\[ d_2^{-1}u_3d_2 u_3= u_3 d_2^{-1}u_3d_2, 
d_3^{-1}u_1d_3 u_1= u_1 d_3^{-1}u_1d_3, 
d_1^{-1}u_2d_1 u_2= u_2 d_1^{-1}u_2d_1\]
\end{proposition}
\begin{proof}
By convention, $ab$ means that we first apply $b$ and then $a$.  
Remark that $u_i=\sigma_{0 p_i(1)}$ is the braid twisting $0$ and $p_i(1)$. 
Moreover,  one proves by recurrence that
$\sigma_{p_i(k)p_i(k+1)}=d_{i-1}^{-k}u_i d_{i-1}^{k}$. Thus 
$AS_{\delta}(Y^*)$ is generated by the $d_j$ since their images generate 
$\Z^2$ and the $d_j$ generate $B_{\infty}$. The presentation for 
$B_{\infty}$  given by Sergiescu (see also \cite{BKL}) yields the 
vertex and edge relations relations 
\[ u_3u_2u_1u_3=u_1u_3u_2u_1=u_2u_1u_3u_2, 
{\rm where }~ u_i=d_{i+1}d_id_{i-1}\]
\[ u_1u_2u_1=u_2u_1u_2, u_3u_2u_3=u_2u_3u_2, u_1u_3u_1=u_3u_1u_3\]
The remaining relations in $B_{\infty}$ are commuting relations between twists
of disjoint supports. The twists on pairs of points at distance one yield 
the relations 
\[ d_2^{-1}u_3d_2 u_1= u_1 d_2^{-1}u_3d_2, 
d_3^{-1}u_1d_3 u_2= u_2 d_3^{-1}u_1d_3, 
d_1^{-1}u_2d_1 u_3= u_3 d_1^{-1}u_2d_1\]
Further $d_i$ commutes with $\sigma_{p_i(1)p_i(2)}$ since they have disjoint 
supports, which lead us to the 
last relations above. 
The interesting phenomenon is that these relations 
actually are sufficient to imply all commutativity relations 
(between arbitrary 
disjoint twists). The proof is a direct calculation (similar to that in 
\cite{Dy}) which we will omit. 

\vspace{0.2cm}
\noindent
Thus the presentation above shows that the subgroup generated by the $u_i$'s
and their conjugates by  the $d_j$ is $B_{\infty}$. Further this subgroup 
$B_{\infty}$ is normal inside  the group given by the presentation above. 
We can provide an infinite presentation of $AS_{\partial}(Y^*)$ 
by Hall's lemma which 
starts from the Sergiescu presentation of $B_{\infty}$ and one of 
$\Z^2$. Then all relations of this infinite presentation are consequences 
of those from the statement. The claim follows.   
\end{proof}
\begin{remark}
It follows from above that 
$B_{\infty}=[AS_{\partial}(Y^*),AS_{\partial}(Y^*)]$ is the commutator 
subgroup since the abelianization of $AS_{\partial}(Y^*)$ is $\Z^2$, as 
in the case studied by Dynnikov. 
\end{remark}

\subsection{$T^*$ is group of homeomorphisms of $S^1$}
 
\begin{proposition}
Let the symbol $^{\diamondsuit}$ denote either $^*$ or $^{\star}$. There exists an embedding 
\[ T^{\diamondsuit}\to {\rm Homeo}_+(S^1)\]
so that $T^{\diamondsuit}$ acts faithfully on the circle. 
\end{proposition} 
\begin{proof}
It is known (see \cite{GS}) that $T$ can be embedded as the subgroup
of piecewise linear homeomorphisms of $S^1=\R/\Z$ which have dyadic 
break points and dyadic derivatives (where defined). This implies that 
$T$ admits a circular order. Furthermore, the group $B_{\infty}$, as all 
finite type braid groups, is left orderable 
(see e.g. \cite{DDRW}, Prop 9.2.7).
By using the exact sequence
\[ 1\to B_{\infty}\to T^{\diamondsuit}\to T\to 1\]
we define a circular order on $T^{\diamondsuit}$, as follows. Let $\pi:T^{\diamondsuit}\to T$
the projection and $x,y,z$ be three elements of $T^{\diamondsuit}$.   
\begin{enumerate}
\item If $\pi(x), \pi(y), \pi(z)$ are distinct then their order in $T^*$
is that of their images in $T$. 
\item If $\pi(x)=\pi(y)\neq \pi(z)$, then $x^{-1}y\in B_{\infty}$ 
which is left orderable. If $x^{-1}y>1$ then $x,y,z$ are positively 
oriented, otherwise it is negatively oriented. 
\item If $\pi(x)=\pi(y)=\pi(z)$ then $x^{-1}y, x^{-1}z\in B_{\infty}$. 
Assume that $1, x^{-1}y, x^{-1}z$ are totally ordered using the 
order in $B_{\infty}$. Then the corresponding $x, y,z$ are 
positively oriented in $T^{\diamondsuit}$.  
\end{enumerate}
This yields a circular order on $T^{\diamondsuit}$ and thus 
there exists an embedding $T^{\diamondsuit}\hookrightarrow   {\rm Homeo}_+(S^1)$. From (\cite{Ca}, Thm.2.2.15) there is a  faithful 
$T^{\diamondsuit}$-action on $S^1$. 
See \cite{Ca} for more details about circular orders and 
related questions.  
\end{proof}

\noindent Any circularly ordered group  $G$ has an embedding 
$G\to {\rm Homeo}_+(S^1)$. A construction due to Ghys and Thurston 
yields a bounded cocycle 
$e$ on $G$ which is the pull-back of the Euler cocycle 
on ${\rm Homeo}_+(S^1)$, and which takes only the values 0 and 1. 
This defines an Euler class $[e]\in H^2(G)$, which is the Euler class 
of the circular order. It is known that $[e]=0$ only if 
$G$ is left ordered. Moreover, in the case of $T^{\diamondsuit}$ with its 
circular ordered defined above, one knows that $T^{\diamondsuit}$ cannot be 
left ordered since it has torsion. This proves that:  
\begin{proposition}
The Euler class 
$[e_{T^{\diamondsuit}}]\in H^2(T^{\diamondsuit})$ is a nontrivial bounded class, whose restriction to $B_{\infty}$ is trivial. 
\end{proposition} 

\noindent Moreover, this Euler class could be also described in 
cohomological terms. Set $e_T$ for the Euler cocycle on the group $T$, namely 
the cocycle inherited from its natural 
embedding into ${\rm Homeo}_+(S^1)$. 

\begin{proposition}
The class $[e_{T^{\diamondsuit}}]\in H^2(T^{\diamondsuit},\Q)$ 
is the pull-back $\pi^*[e_T]$ by the projection $\pi:T^{\diamondsuit}\to T$. 
\end{proposition}
\begin{proof}
According to a result due to E.Ghys, S.Jekel, and W.Thurston we have 
\[ [c] = 2 [e]  \]
where $c$ is the order cocycle  defined by Thurston (see \cite{Ca}, 
Construction 2.3.4) 
and $e$ is the Euler cocycle of 
a circularly ordered group.  If the group $G$ acts faithfully on $S^1$ let 
us choose a point $a\in S^1$ with trivial stabilizer. 
Recall that $c$ is defined (as a homogeneous cocycle) by means of 
\[ c(g_0:g_1:g_2)=\left\{\begin{array}{cccl}
 1  & {\rm if} & (g_0(a),g_1(a),g_2(a)) & {\rm is ~positively ~oriented}\\
 -1 & {\rm if} & (g_0(a),g_1(a),g_2(a)) & {\rm is ~negatively ~ oriented}\\
0   & {\rm if} & (g_0(a),g_1(a),g_2(a)) & {\rm is ~ degenerate}
\end{array}\right.
\]

\vspace{0.1cm}
\noindent
We claim that
$p^*([c_T])=[c_{T^{\diamondsuit}}]$. 
More generally, if $p:G\to H$ is a  monotone homomorphism 
with  left orderable kernel between circularly orderable groups 
$G$ and $H$ then $p^*[c_H] = [c_G]$. 

\vspace{0.1cm}
\noindent
In fact, we have 
\[ c_G-p^*c_H = \partial w \]
where $w$ is the following 1-cocycle (in homogeneous
coordinates):

\[w(g_0:g_1)=\left\{\begin{array}{ccccc}
 0  & {\rm if} & p(g_0)\neq p(g_1) &           &\\
 1  & {\rm if} & p(g_0)=p(g_1)    & {\rm and} & g_0^{-1}g_1 <1\\
-1  & {\rm if} & p(g_0)=p(g_1)    & {\rm and} & g_0^{-1}g_1 >1 \\
 0  & {\rm if} &  g_0=g_1          &            &
\end{array}\right. \]
This implies that $p^*[e_H]=[e_G]$ up to 2-torsion, as claimed.  
\end{proof}

 \section{The complex $\C$} 
 The following of the article is devoted to the proof that $\T$ is finitely generated, by constructing a simply connected 
cellular complex $\CR$ on which $\T$ acts cocompactly. Each orbit of 2-cells of this complex will thus correspond to a relation in $\T$. This
will enable us to provide an explicit presentation for $\T$. We first introduce an auxiliary complex $\C$, whose
simple connectivity is not too difficult to prove, but which is not finite modulo $\T$. The complex $\CR$ will be
a subcomplex of $\C$.
We shall
finally explain how an analogous construction applies to $T^{\star}$. 

 \subsection{Vertices and Edges}

\begin{figure}
\begin{center}
\begin{picture}(0,0)%
\includegraphics{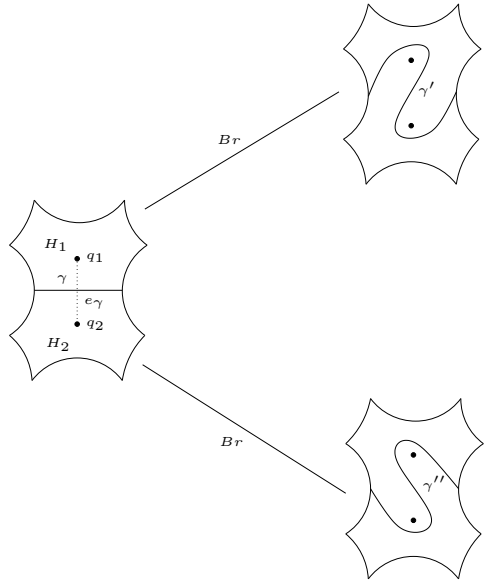}%
\end{picture}%
\setlength{\unitlength}{1243sp}%
\begingroup\makeatletter\ifx\SetFigFont\undefined%
\gdef\SetFigFont#1#2#3#4#5{%
  \reset@font\fontsize{#1}{#2pt}%
  \fontfamily{#3}\fontseries{#4}\fontshape{#5}%
  \selectfont}%
\fi\endgroup%
\begin{picture}(9466,11491)(1163,-10959)
\put(1846,-4336){\makebox(0,0)[lb]{\smash{{\SetFigFont{5}{6.0}{\rmdefault}{\mddefault}{\updefault}{$H_1$}%
}}}}
\put(1891,-6316){\makebox(0,0)[lb]{\smash{{\SetFigFont{5}{6.0}{\rmdefault}{\mddefault}{\updefault}{$H_2$}%
}}}}
\put(5311,-2266){\makebox(0,0)[lb]{\smash{{\SetFigFont{5}{6.0}{\rmdefault}{\mddefault}{\updefault}{$Br$}%
}}}}
\put(5356,-8296){\makebox(0,0)[lb]{\smash{{\SetFigFont{5}{6.0}{\rmdefault}{\mddefault}{\updefault}{$Br$}%
}}}}
\put(9316,-1276){\makebox(0,0)[lb]{\smash{{\SetFigFont{5}{6.0}{\rmdefault}{\mddefault}{\updefault}{$\gamma'$}%
}}}}
\put(2656,-5461){\makebox(0,0)[lb]{\smash{{\SetFigFont{5}{6.0}{\rmdefault}{\mddefault}{\updefault}{$e_{\gamma}$}%
}}}}
\put(2116,-5011){\makebox(0,0)[lb]{\smash{{\SetFigFont{5}{6.0}{\rmdefault}{\mddefault}{\updefault}{$\gamma$}%
}}}}
\put(2701,-5911){\makebox(0,0)[lb]{\smash{{\SetFigFont{5}{6.0}{\rmdefault}{\mddefault}{\updefault}{$q_2$}%
}}}}
\put(9406,-9151){\makebox(0,0)[lb]{\smash{{\SetFigFont{5}{6.0}{\rmdefault}{\mddefault}{\updefault}{$\gamma''$}%
}}}}
\put(2701,-4561){\makebox(0,0)[lb]{\smash{{\SetFigFont{5}{6.0}{\rmdefault}{\mddefault}{\updefault}{$q_1$}%
}}}}
\end{picture}%

\caption{definition of the $Br$-move}\label{Br}
\end{center}
\end{figure}

 The complex $\C$ is a 2-dimensional cellular complex, whose vertices
 correspond to the elements of $\mathfrak{R}^*$ (cf. Notation \ref{not}). The (unoriented) edges correspond to
 ``moves'' of two types:

 \begin{enumerate}
 \item Associativity move, or $A$-move:

\begin{figure}
\begin{center}
\begin{picture}(0,0)%
\includegraphics{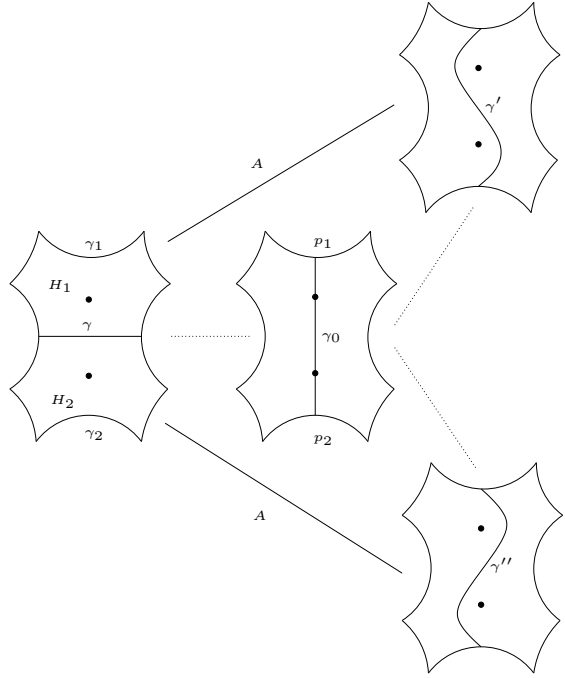}%
\end{picture}%
\setlength{\unitlength}{1450sp}%
\begingroup\makeatletter\ifx\SetFigFont\undefined%
\gdef\SetFigFont#1#2#3#4#5{%
  \reset@font\fontsize{#1}{#2pt}%
  \fontfamily{#3}\fontseries{#4}\fontshape{#5}%
  \selectfont}%
\fi\endgroup%
\begin{picture}(9466,11491)(1163,-10959)
\put(9316,-1276){\makebox(0,0)[lb]{\smash{\SetFigFont{5}{6.0}{\rmdefault}{\mddefault}{\updefault}{$\gamma'$}%
}}}
\put(1846,-4336){\makebox(0,0)[lb]{\smash{\SetFigFont{5}{6.0}{\rmdefault}{\mddefault}{\updefault}{$H_1$}%
}}}
\put(1891,-6316){\makebox(0,0)[lb]{\smash{\SetFigFont{5}{6.0}{\rmdefault}{\mddefault}{\updefault}{$H_2$}%
}}}
\put(5356,-8296){\makebox(0,0)[lb]{\smash{\SetFigFont{5}{6.0}{\rmdefault}{\mddefault}{\updefault}{$A$}%
}}}
\put(2476,-6856){\makebox(0,0)[lb]{\smash{\SetFigFont{5}{6.0}{\rmdefault}{\mddefault}{\updefault}{$\gamma_2$}%
}}}
\put(5311,-2266){\makebox(0,0)[lb]{\smash{\SetFigFont{5}{6.0}{\rmdefault}{\mddefault}{\updefault}{$A$}%
}}}
\put(6391,-3571){\makebox(0,0)[lb]{\smash{\SetFigFont{5}{6.0}{\rmdefault}{\mddefault}{\updefault}{$p_1$}%
}}}
\put(6391,-6946){\makebox(0,0)[lb]{\smash{\SetFigFont{5}{6.0}{\rmdefault}{\mddefault}{\updefault}{$p_2$}%
}}}
\put(6526,-5191){\makebox(0,0)[lb]{\smash{\SetFigFont{5}{6.0}{\rmdefault}{\mddefault}{\updefault}{$\gamma_0$}%
}}}
\put(2476,-3616){\makebox(0,0)[lb]{\smash{\SetFigFont{5}{6.0}{\rmdefault}{\mddefault}{\updefault}{$\gamma_1$}%
}}}
\put(2431,-5011){\makebox(0,0)[lb]{\smash{\SetFigFont{5}{6.0}{\rmdefault}{\mddefault}{\updefault}{$\gamma$}%
}}}
\put(9406,-9151){\makebox(0,0)[lb]{\smash{\SetFigFont{5}{6.0}{\rmdefault}{\mddefault}{\updefault}{$\gamma''$}%
}}}
\end{picture}

\caption{Definition of the A-move}\label{A}
\end{center}
\end{figure}

\vspace{0.1cm} 
\noindent 
 Let $\mathfrak{r}$ be a rigid structure, and
 $\gamma$ an arc of  $\mathfrak{r}$ which separates two hexagons, say $H_1$ and $H_2$. Let  $\gamma_1$
 (resp. $\gamma_2$) be the side of  $H_1$ (resp. $H_2$) contained in the
boundary
 of $D^*$, and disjoint from $\gamma$. Choose two points $p_1\in \gamma_1$ and
 $p_2\in\gamma_2$. Define the simple arc $\gamma_0$, which first connects $p_1$ to
the
 puncture of $H_1$ (remaining inside $H_1$), next connects the puncture of
 $H_1$ to the puncture of  $H_2$, crossing once and transversely the arc $\gamma$, and
 finally connects the puncture of  $H_2$ to $p_2$ (remaining inside $H_2$).
Define an arc  by deforming $\gamma_0$ around the two
punctures, in such a way that it avoids them and separates them. There
are exactly two isotopy classes of such arcs, say $\gamma'$ and $\gamma''$, with free extremities $p_1$
and $p_2$ (see Figure \ref{A}).

\vspace{0.1cm} 
\noindent 
One says that the rigid structures  ${\mathfrak r}'$ and ${\mathfrak r}''$, obtained from
${\mathfrak r}$ by only changing $\gamma$ into $\gamma'$ or $\gamma''$, respectively, are obtained from
${\mathfrak r}$ by an
{\it $A$-move on $\gamma$}. A pair as $\{\gamma,\gamma'\}$ or $\{\gamma,\gamma''\}$ determines an {\it edge of type $A$} of $\C$. Note that there exist exactly two $A$-moves on $\gamma$.

\vspace{0.1cm}\noindent 
 \item Braiding move, or $Br$-move:

\vspace{0.2cm}\noindent 
Let $\mathfrak{r}$ be a rigid structure, and
 $\gamma$ an arc of $\mathfrak{r}$ which separates two hexagons, say $H_1$ and $H_2$. Let
 $e_{\gamma}$ be a simple arc connecting $q_1$ to $q_2$ (the punctures of $H_1$ and $H_2$, respectively), crossing $\gamma$
  once and transversely, 
 and contained in $H_1\cup H_2$. Such an arc is uniquely defined, up to isotopy (fixing $q_1$ et $q_2$). Let $\sigma$ be the
 positive braid determined by $e_{\gamma}$, which permutes $q_1$ and
 $q_2$ (cf. Definition \ref{defbraid}). Let
 $\gamma'=\sigma(\gamma)$ be the image of $\gamma$ by $\sigma$, and $\gamma''=\sigma^{-1}(\gamma)$ (see Figure \ref{Br}).

\vspace{0.1cm}\noindent 
One says that the rigid structures ${\mathfrak r}'$ and ${\mathfrak r}''$, obtained from ${\mathfrak r}$ by
 only changing $\gamma$ into $\gamma'$ or $\gamma''$, respectively, are obtained from
${\mathfrak r}$ by a {\it $Br$-move on $\gamma$}. A pair as $\{\gamma,\gamma'\}$ or $\{\gamma,\gamma''\}$ determines
an {\it edge of type $Br$} of $\C$. Note that there exist exactly two $Br$-moves on $\gamma$.
\end{enumerate}

\begin{definition}
The tree ${\cal T}_{\mathfrak r}$ of a rigid structure $\mathfrak r$ of $D^*$ is the planar 
tree whose vertices are the punctures of $D^*$ and whose edges are the arcs $e_{\gamma}$ as above, for every arc $\gamma$ of
 $\mathfrak r$.
\end{definition}

\vspace{0.1cm}\noindent 
Note that ${\cal T}_{{\mathfrak r}^*}$ is the tree of $D^*$.

\begin{remark}[Orientation of the edges of type $Br$] Say that an edge $\gamma\stackrel{Br}{\longrightarrow} \gamma'$ is positively oriented if
 $\gamma'$ turns on the left when it approaches the arc $e_{\gamma}$. This means that the braiding $\sigma$ (on the punctures $q_1$ and $q_2$) such that 
 $\gamma'=\sigma(\gamma)$, is  positive. From now on, a positive braiding will always be denoted
  by a letter without negative exponent, such as $\sigma$, $\sigma_1$, etc., while $\sigma^{-1}$, $\sigma_1^{-1}$, etc.,  will refer to negative braidings. On
   Figure \ref{Br}, the edge $\gamma\stackrel{Br}{\rightarrow}\gamma'$ is positively oriented, while the edge  $\gamma\stackrel{Br}{\rightarrow}\gamma''$ 
   is negatively oriented.

\vspace{0.1cm}\noindent 
   On the contrary, there is no canonical orientation for the edges of type $A$.
\end{remark}

\subsection{2-cells}\label{2-cells}

The 2-cells of $\C$ are of the following types:

\vspace{0.2cm}\noindent 
1. Cells $AA=Br$.

\begin{figure}
\begin{center}

\includegraphics{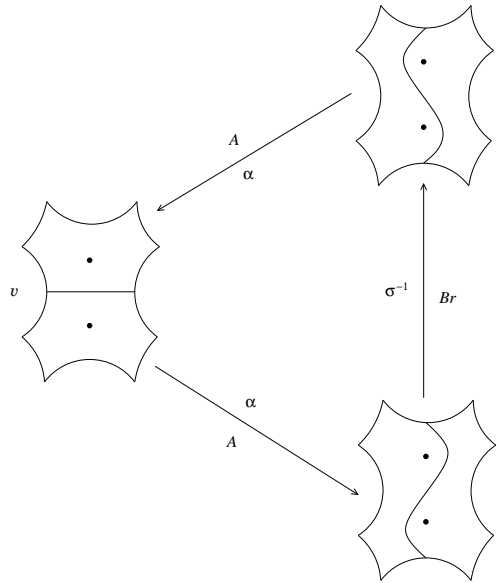}

\caption{Cycle $AA=Br$}\label{AA=Br}
\end{center}
\end{figure}

\vspace{0.2cm}\noindent 
Let ${\mathfrak r}$ be a rigid structure, and $\gamma$ an arc of  ${\mathfrak r}$. As we have seen above,
there are two edges of type $A$, connecting ${\mathfrak r}$ to ${\mathfrak
  r}'$ and  ${\mathfrak r}''$. The vertices ${\mathfrak  r}'$ and ${\mathfrak
  r}''$ are connected by an edge of type $Br$. Thus, one fills in the cycle of
three edges $A,A$ and $Br$, by a 2-cell, which is said of {\it type $AA=Br$} (see Figure \ref{AA=Br}).

\vspace{0.2cm}\noindent 
2. Cells of commutation of $A$-moves, $A_1A_2=A_2A_1$.

\begin{figure}
\begin{center}
\includegraphics{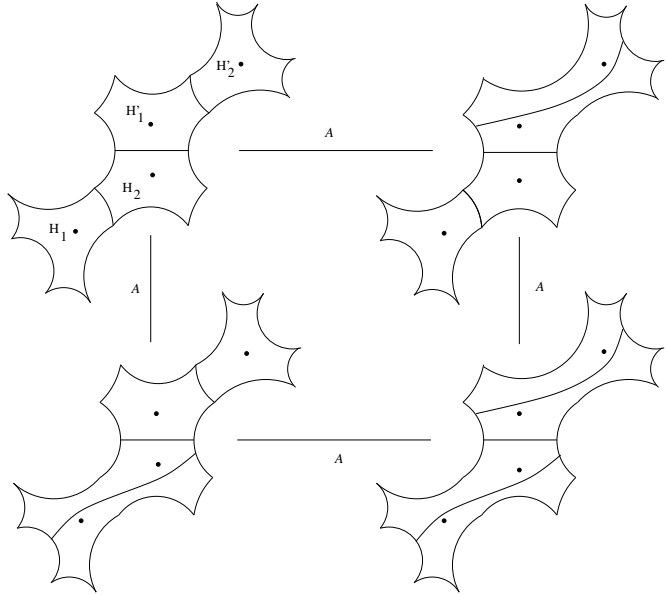}
\caption{Cycle of type $DC_1$}\label{DC1}
\end{center}
 \end{figure}

\vspace{0.2cm}\noindent 
Let $\mathfrak{r}$ be a rigid structure, and
 $H_{1}, H_{2},H'_{1}$ and  $H'_{2}$ be four distinct hexagons of $\mathfrak{r}$, such that $H_1$ and $H_2$ (resp. $H'_1$ and $H'_2$) share a common side 
 $\gamma$ (resp. $\gamma'$). The commutation of the two $A$-moves, along $\gamma$ and $\gamma'$, respectively, generates a square cycle. The point
  (to be elucidated in Proposition \ref{HT}) is that we only need to fill in the squares of two special types: 

\vspace{0.2cm} 
$\bullet$ Cells $DC_1$: $H_2$ and $H_1'$ share a common side, see Figure \ref{DC1}.
    
\begin{figure}
\begin{center}
\includegraphics{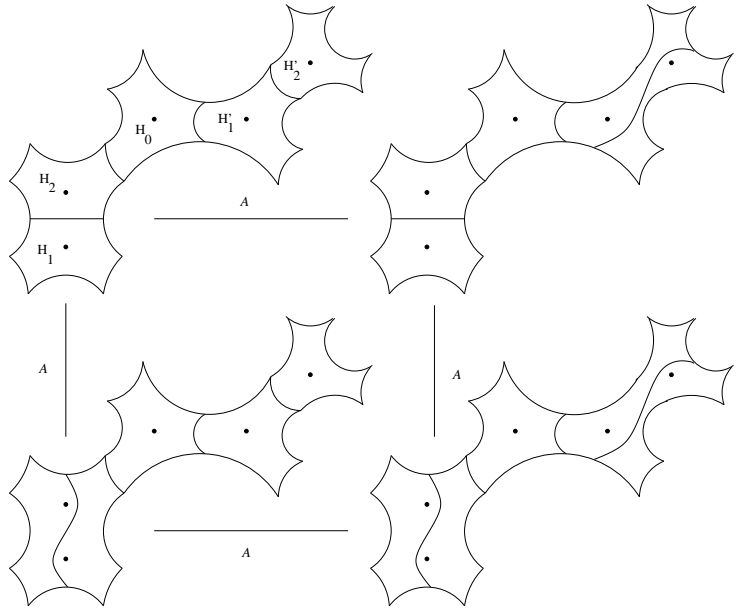}
\caption{Cycle of type $DC_2$}\label{DC2}
\end{center}
\end{figure}

\vspace{0.2cm}
$\bullet$ Cells $DC_2$: $H_2$ and $H_1'$ are separated by a hexagon $H_{0}$, see Figure \ref{DC2}.

\vspace{0.2cm}\noindent 
3. Pentagonal cells.

\begin{figure}
\begin{center}
\includegraphics{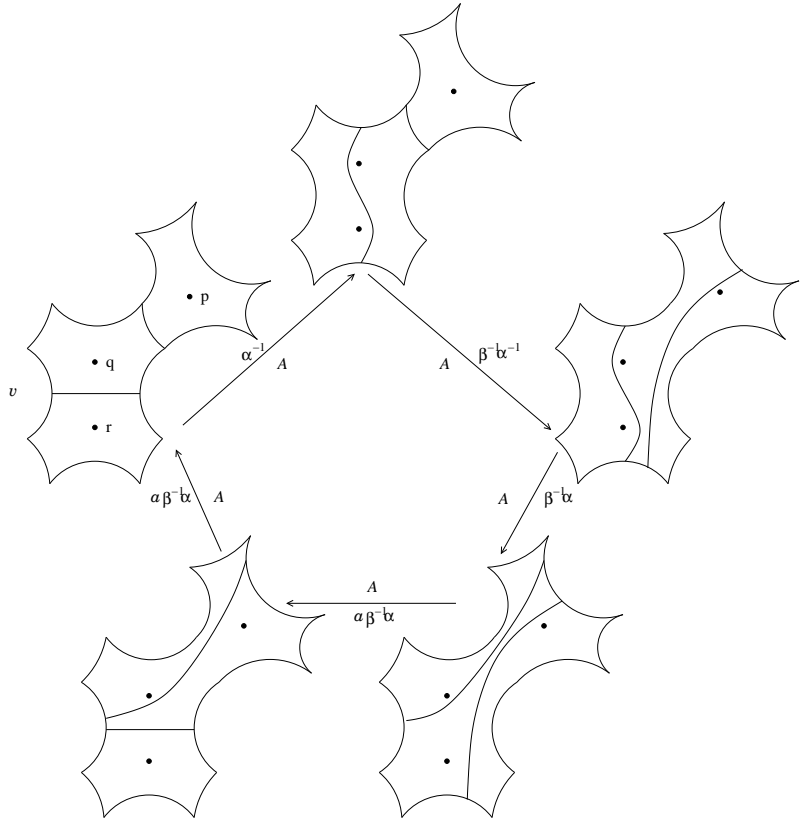}
\caption{Pentagonal cycle}\label{penta}
\end{center}
\end{figure}

\vspace{0.2cm}\noindent 
Let $\mathfrak{r}$ be a rigid structure, and $H_1$, $H_2$ and $H_3$ be 
three hexagons of $\mathfrak{r}$, such that  $H_1$ and $H_2$ are adjacent along a side $\gamma$, and $H_3$ 
and $H_2$ are adjacent along a side $\delta$. There is a cycle of five $A$-moves, which only involves the arcs $\gamma$
 and $\delta$ of $\mathfrak{r}$, see Figure \ref{penta}. It is filled in, producing a 2-cell of {\it pentagonal type}.

\vspace{0.2cm}\noindent 
4. Cells coming from the presentation of the braid group.

\vspace{0.2cm}\noindent 
Recall that there is a general theorem of V. Sergiescu (\cite{Se}) which can be used to provide a presentation for $B_{\infty}$
 with  generators the positive braidings 
along the edges of the tree of $D^*$ or the tree ${\cal T}_{\mathfrak r}$ of any rigid structure ${\mathfrak r}$. Let $\mathfrak{r}$ be a rigid structure, and ${\cal T}_{\mathfrak r}$ be its tree.
  
\vspace{0.2cm}
a) Hexagonal cells. Let $e_1$ and $e_2$ be two edges of ${\cal T}_{\mathfrak r}$, which are incident to a puncture $p$. Let $\sigma_1$ and 
$\sigma_2$ be the braidings along $e_1$ and $e_2$, respectively. Then $\sigma_1\sigma_2  \sigma_1=\sigma_2\sigma_1\sigma_2$.
 Let $Br_1$ and $Br_2$ denote the braiding moves, corresponding to $\sigma_1$ and $\sigma_2$, respectively. In terms of $Br$-moves, the braid relation becomes
$$Br_2Br_1Br_2=Br_1Br_2Br_1,$$
and one adds a 2-cell to fill in the cycle of those 6 braiding moves.

\vspace{0.2cm} 
b) Octagonal cells. Let $e_1$, $e_2$ and $e_3$ be the three edges which are incident to a puncture $p$. Suppose that their
 enumeration respects the cyclic counterclockwise orientation
 of the planar surface around $p$. Using notations as in a), one has the  relation  $\sigma_1\sigma_2  \sigma_3\sigma_1=\sigma_2\sigma_3\sigma_1\sigma_{2}= \sigma_3\sigma_1\sigma_2\sigma_{3}$. In terms of $Br$-moves,
 this gives 
$$Br_3Br_2Br_1Br_3= Br_2Br_1Br_3Br_2= Br_1 Br_3Br_2Br_1,$$
and one adds 2-cells to fill in the corresponding cycles of 8 braidings.

\vspace{0.2cm} 
c) Squares. Let $e_1$ and $e_2$ two disjoint edges. Then $\sigma_1\sigma_2=\sigma_2\sigma_1$. In terms of braiding moves, this gives 
$$Br_1Br_2=Br_2Br_1,$$
and one adds a 2-cell to fill in this square cycle.

\begin{figure}
\begin{center}
\includegraphics{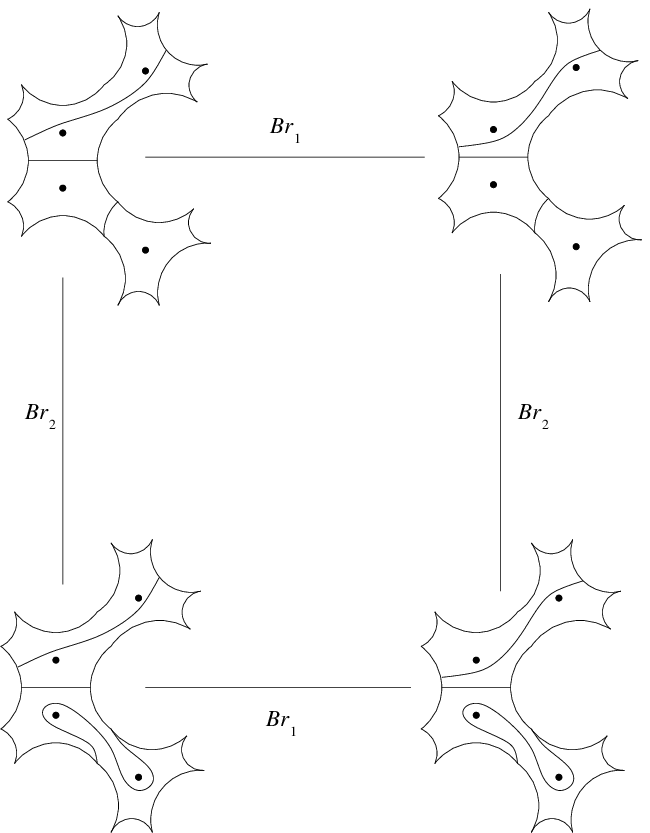}
\caption{Commutation $Br_1Br_2=Br_2Br_1$ of level 6}\label{Brlev6}
\end{center}
 \end{figure}

\begin{figure}
\begin{center}
\includegraphics{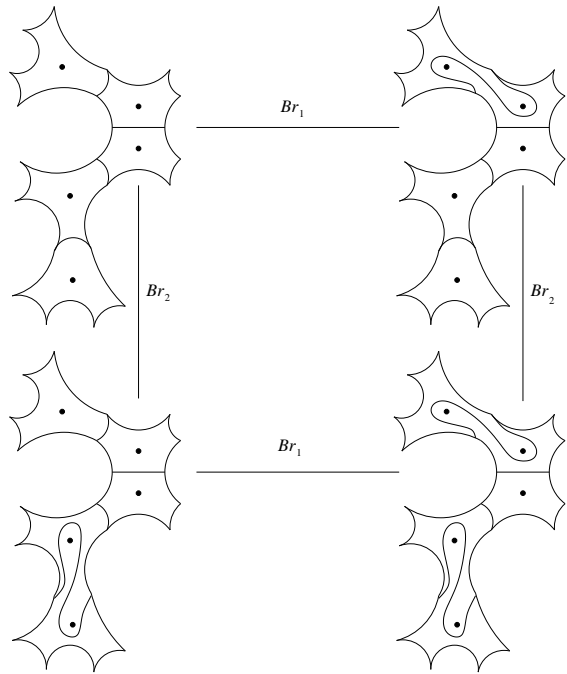}
\caption{First commutation $Br_{1}Br_{2}=Br_{2}Br_{1}$ of level 7}\label{Brlev7}
\end{center}
\end{figure}
 
\begin{figure}
\begin{center}
\includegraphics{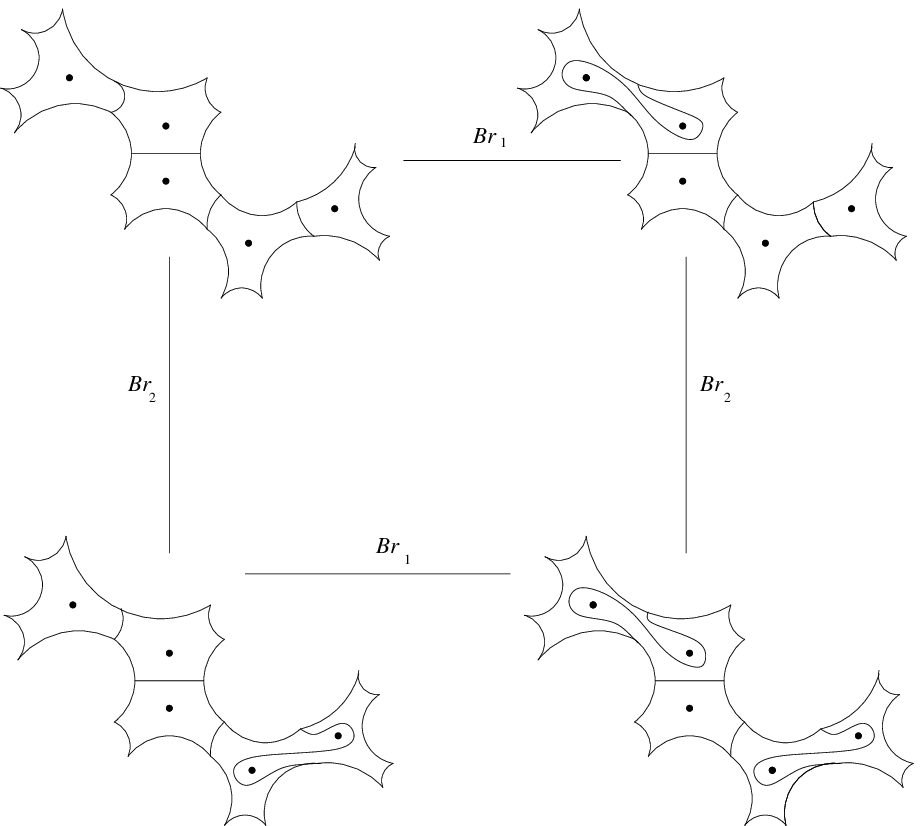}
\caption{Second commutation $Br_{1}Br_{2}=Br_{2}Br_{1}$ of level 7}\label{Brlev7bis}
\end{center}
\end{figure}  

\vspace{0.2cm}\noindent 
5. Cells of commutation of $A$-moves with $Br$-moves.
 
\vspace{0.2cm}\noindent 
Let $\mathfrak{r}$ be a rigid structure. An $A$-move along an arc $\gamma$ commutes with a $Br$-move along an edge $e$ of ${\cal T}_{\mathfrak r}$  if $\gamma$ and $e$ are disjoint. 
Thus, there is a square cycle of the form $$ABr=BrA$$
which one fills in by a 2-cell.

\begin{figure}
\begin{center}
\includegraphics{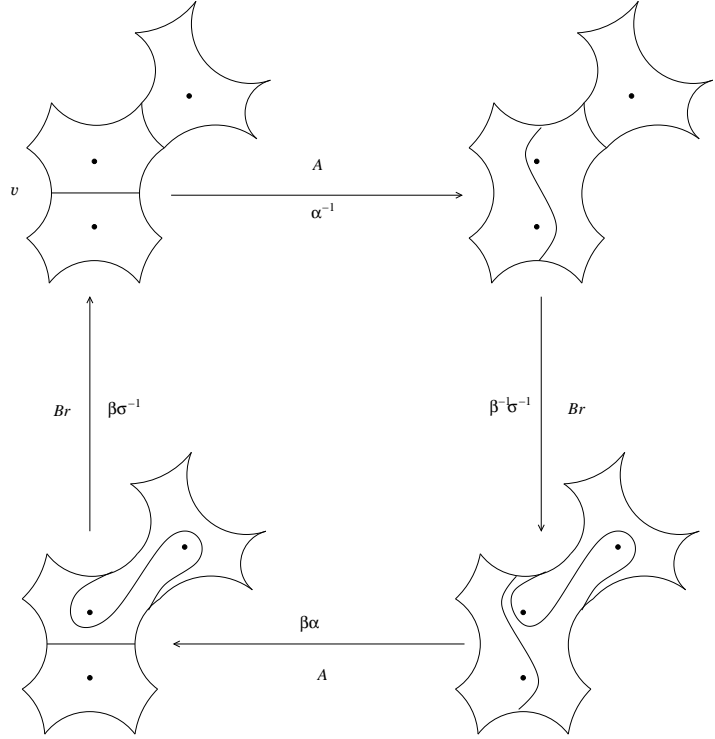}
\caption{Cycle $BrA=ABr$ of level 5}\label{BrA=ABr5}
\end{center}
\end{figure}

\vspace{0.2cm}\noindent 
We note that the minimum level (see the definition below) for such a cell is 5, see Figure \ref{BrA=ABr5}. 

\begin{definition}
Let $\omega$ be a 2-cell, and $r$ a vertex of the boundary cycle $\partial \omega$. The vertices of  $\partial \omega$ differ from $r$ 
by a finite number of arcs  $\gamma$. The support of $\omega$ is the minimal connected 
subsurface of $D^*$ which is a union of hexagons of $r$ and contains all the arcs $\gamma$. The level of the 2-cell
$\omega$ is the number of arcs of $r$ which belong to the boundary of the support of $\omega$.
\end{definition}

\vspace{0.2cm}\noindent 
The description of $\C$ is now complete, and the following is obvious:

\begin{proposition}
The complex $\C$ is a $\T$-complex. 
\end{proposition}

\subsection{Connectivity of $\C$} 

We first recall  a useful lemma of algebraic topology (\cite{BK2}, prop. 6.2, see also a variant 
of it in \cite{FG}), which we have used already in \cite{FK}.     
         
\begin{lemma}\label{sta}     
Let ${\cal M}$ and ${\cal C}$ be two $CW$-complexes of dimension          
  $2$, with oriented edges, and $f :{\cal M}^{(1)}\rightarrow {\cal C}^{(1)}$          
  be a cellular map between their 1-skeletons, which is surjective on $0$-cells and          
  $1$-cells. Suppose that:            
\begin{enumerate}          
\item ${\cal C}$ is connected and simply connected;          
\item For each vertex $c$ of ${\cal C}$,            
$f^{-1}(c)$ is connected and simply connected;          
\item Let $c_1\stackrel{e}{\longrightarrow } c_2$ be an oriented edge 
  of ${\cal C}$,          
  and let $m_1 '\stackrel{e'}{\longrightarrow } m_2 '$ and $m_1''     
  \stackrel{e''}{\longrightarrow } m_2''$ be two lifts in ${\cal M}$. Then          
  we can find two paths $m_1 '\stackrel{p_1}{\longrightarrow } m_1 ''$ in          
  $f^{-1}(c_1)$ and  $m_2 '\stackrel{p_2}{\longrightarrow } m_2 ''$ in          
  $f^{-1}(c_2)$ such that the loop        
          
\begin{center}      
\begin{picture}(0,0)%
\includegraphics{dia.pstex}%
\end{picture}%
\setlength{\unitlength}{2486sp}%
\begingroup\makeatletter\ifx\SetFigFont\undefined%
\gdef\SetFigFont#1#2#3#4#5{%
  \reset@font\fontsize{#1}{#2pt}%
  \fontfamily{#3}\fontseries{#4}\fontshape{#5}%
  \selectfont}%
\fi\endgroup%
\begin{picture}(1845,1650)(1396,-2041)     
\put(2926,-2041){\makebox(0,0)[lb]{\smash{\SetFigFont{8}{9.6}{\rmdefault}{\mddefault}{\updefault}$m_2''$}}}     
\put(1576,-2041){\makebox(0,0)[lb]{\smash{\SetFigFont{8}{9.6}{\rmdefault}{\mddefault}{\updefault}$m_1''$}}}     
\put(1396,-1366){\makebox(0,0)[lb]{\smash{\SetFigFont{8}{9.6}{\rmdefault}{\mddefault}{\updefault}$p_1$}}}     
\put(3241,-1366){\makebox(0,0)[lb]{\smash{\SetFigFont{8}{9.6}{\rmdefault}{\mddefault}{\updefault}$p_2$}}}     
\put(2431,-556){\makebox(0,0)[lb]{\smash{\SetFigFont{8}{9.6}{\rmdefault}{\mddefault}{\updefault}$e'$}}}     
\put(2431,-1861){\makebox(0,0)[lb]{\smash{\SetFigFont{8}{9.6}{\rmdefault}{\mddefault}{\updefault}$e''$}}}     
\put(2971,-736){\makebox(0,0)[lb]{\smash{\SetFigFont{8}{9.6}{\rmdefault}{\mddefault}{\updefault}$m_2'$}}}     
\put(1621,-736){\makebox(0,0)[lb]{\smash{\SetFigFont{8}{9.6}{\rmdefault}{\mddefault}{\updefault}$m_1'$}}}     
\end{picture}     
       
\end{center}            
is contractible in ${\cal M}$;          
\item For any $2$-cell $X$ of ${\cal C}$, its boundary $\partial X$ can be          
 lifted to a contractible loop of ${\cal M}$.          
\end{enumerate}          
Then ${\cal M}$ is connected and simply connected.         
\end{lemma}  

\vspace{0.2cm}\noindent 
Recall that $D$ is the surface $D^*$ viewed without its punctures. We will use the Lemma to study a certain cellular
 map $f:\C\rightarrow {\cal HT}_{red}(D).$

\begin{definition}
 The reduced Hatcher-Thurston complex ${\cal HT}_{red}(D)$ is a 2-dimensional cellular complex
 whose vertices are the rigid structures of $D$, whose edges correspond to $A$-moves, and whose 2-cells are of three types:
  $DC_1$, $DC_2$, and pentagonal cells.  The definition of the $A$-move in ${\cal HT}_{red}(D)$ is deduced from the
  definition of the $A$-move in $\C$ by forgetting the punctures. 
\end{definition}

\vspace{0.2cm}\noindent 
\noindent Note that, if $\gamma$ is an arc of a rigid structure of $D$, there is a unique $A$-move on $\gamma$.

\begin{remark}\label{rem}
\begin{enumerate}   
\item If $\Sigma_{0,\infty}$ is the surface without boundary obtained by gluing along their boundaries two copies of $D$ with opposite 
orientations, 
then  ${\cal HT}_{red}(D)$ is a subcomplex of the {\it reduced Hatcher-Thurston complex} ${\cal HT}_{red}(\Sigma_{0,\infty})$ of the surface 
 $\Sigma_{0,\infty}$, 
as it appears in \cite{FK}, Definition 5.2. The argument used in \cite{FK} to prove that ${\cal HT}_{red}(\Sigma_{0,\infty})$ is connected and simply connected 
actually reduces to proving that ${\cal HT}_{red}(D)$ is connected and simply connected. The point of that
  proof (see Proposition 5.5 in \cite{FK}) is that there is a surjection of the Cayley complex
   of Thompson's group $T$, for the presentation with generators $\alpha$ and $\beta$, onto
    the complex ${\cal HT}_{red}(D)$. This is essentially used to show that the square 
    cycles generated by the commutations of any two A-moves are filled in by 2-cells of types $DC_1$, $DC_2$, and by pentagons.
\item The complex ${\cal HT}_{red}(D)$ is a $T$-complex, and ${\cal HT}_{red}(D)/T$ has one vertex, one edge, and three 2-cells: the two squares $DC_1$ and $DC_2$, and the pentagon.
\end{enumerate}
 
\end{remark}

\vspace{0.2cm}\noindent 
The following is obvious:
\begin{proposition}
There is a well defined cellular map
$$f:\C\rightarrow {\cal HT}_{red}(D),$$
which is induced by forgetting the punctures. The map $f$ is $(\T, T)$-equivariant.
\end{proposition} 

\begin{definition}
The $\T$-type of a 2-cell $\omega$ of $\C$ is its image in $\C/\T$. The $T$-type of $\omega$ is $f(\omega)$ mod $T$ in ${\cal HT}_{red}(D)/T$. 
\end{definition}

\begin{proposition}\label{types}
There is exactly one $T$-type of 2-cells $DC_1$, one $T$-type of 2-cells $DC_2$, and one $T$-type of pentagonal 2-cells. 
Each $T$-type of 2-cell corresponds to finitely many different $\T$-types. In other words, if $\omega$ is a 2-cell in ${\cal HT}_{red}(D)$,
 then the set of 2-cells in $\C$ which are the preimages of $\omega$ by $f$ is finite modulo $\T$.
\end{proposition}   
        
\begin{proof}
The first assertion was already mentioned in Remark \ref{rem}, 2. The second is related to the fact that an edge of type $A$ in ${\cal HT}_{red}(D)$ admits two lifts in $\C$ with the same origin (see Figure \ref{AA=Br}). Therefore, a 2-cell in
${\cal HT}_{red}(D)$ bounded by a cycle of $n$ edges admits {\it at most} $2^{n-1}$ lifts in $\C$ based at the same
origin.

\end{proof}

\begin{proposition}\label{HT}
The complex $\C$ is connected and simply connected.
\end{proposition}

\begin{proof}
Let us apply Lemma \ref{sta} to the map $f$. Condition {\it 1} is fulfilled. The preimage by $f$ of a vertex $\mathfrak r$ of
 ${\cal HT}_{red}(D^*)$
 is isomorphic to the Cayley complex of the group $B_{\infty}$, for the presentation of Sergiescu associated to the tree 
 ${\cal T}_{\mathfrak r}$ of ${\mathfrak r}$. Consequently,
  it is connected and simply connected, and condition {\it 2} of Lemma \ref{sta} is fulfilled.

\vspace{0.2cm}\noindent 
   Let us examine condition {\it 3}. The edges $e'$
   and $e''$ are of type $A$. Since $m'_{1}$ (resp. $m'_{2}$) is connected to $m''_{1}$ (resp. $m''_{2}$) by a sequence
of edges of type $Br$, it suffices to consider the case 
   when $p_{1}$ is an edge of type $Br$. But this forces $p_{2}$ either to be trivial (the loop $p_{1}e''e'^{-1}$ bounds a 2-cell $AA=Br$) or 
   to be an edge of type $Br$ (the loop $e'p_{2}e''^{-1}p_{1}^{-1}$ bounds a 2-cell $ABr=BrA$).

\vspace{0.2cm}\noindent 
Condition {\it 4} is obviously fulfilled, by definition of $\C$. To conclude, $\C$ is connected and simply connected.   

\end{proof}

\section{The reduced complex $\CR$ and a presentation for $\T$}

\subsection{Simple connectivity of $\CR$}

\begin{definition}
The reduced complex $\CR$ is the subcomplex of $\C$ which has the same 1-skeleton as $\C$. The 2-cells are of the following types:
\begin{enumerate}
\item $AA=Br$ (Figure \ref{AA=Br}), $DC_1$ of Figure \ref{DC1}, $DC_2$ of Figure \ref{DC2}, pentagon of Figure \ref{penta};
\item $ABr=BrA$ of level 5 of Figure \ref{BrA=ABr5}, level 6 of Figure \ref{ineqABr=BrA}, and level 7 of Figures \ref{BrA=ABr7}
 and \ref{BrA=ABr7bis};
 \item  $Br_{1}Br_{2}=Br_{2}Br_{1}$ coming from the commutation of certain braidings with disjoint supports: 
cells of level 6 of the $T^*$-type of Figure \ref{Brlev6}, and cells of level 7 the $T^*$-type of Figure \ref{Brlev7} and of the $T^*$-type of Figure \ref{Brlev7bis});
 \item cells coming from the braid group: $Br_{1}Br_{2}Br_{2}=Br_{2}Br_{1}Br_{2}$ (subsection \ref{2-cells}, 4.a) and 
 $Br_3Br_2Br_1Br_3= Br_2Br_1Br_3Br_2= Br_1 Br_3Br_2Br_1$ (subsection \ref{2-cells}, 4.b).
 \end{enumerate}

\end{definition}

  \vspace{0.2cm}\noindent  
The point is that $\CR$, contrary to $\C$, contains finitely many $\T$-types of cells $ABr=BrA$ and $Br_{1}Br_{2}=Br_{2}Br_{1}$. 
It follows that the quotient $\CR/\T$ is a finite complex. Moreover:

\begin{proposition}\label{main}
The complex $\CR$ is connected and simply connected.
\end{proposition}

\begin{proof}
Since $\C$ is connected and has the same 1-skeleton as $\CR$, the latter is connected as well. To prove the simple connectivity
 of $\CR$ from the simple connectivity of $\C$, it suffices to check that the cycles bounding the 2-cells which belong to $\C$ 
 but not to $\CR$ may be filled in by some combinations of 2-cells of $\CR$ only.

\vspace{0.2cm}\noindent 
 Note first that for each of the three $T$-types ($DC_{1}$, $DC_{2}$ or pentagon), we have selected a unique $\T$-type of lift in 
 $\CR$ (compare with Proposition \ref{types}). However:
 
 \begin{lemma}\label{le1}
 Let $\omega$ be any 2-cell of $\C$ of $T$-type $DC_{1}$, $DC_{2}$, or of pentagonal $T$-type. Then
  $\partial \omega$ is filled in by 2-cells which belong to $\CR$, hence is homotopically trivial in $\CR$.
\end{lemma} 

\begin{proof}

Let us introduce the following terminology. Suppose that the boundary of a 2-cell $\omega'$ is filled in by some 2-cells $\omega$, $\omega_{1},\ldots,\omega_{n}$. Then we will say that $\omega'$ is equivalent to $\omega$ modulo $\omega_{1},\ldots,\omega_{n}$.

\begin{itemize}
\item  Let us consider the $\T$-types of cells of $T$-type $DC_{1}$. The only $\T$-type which belongs to $\CR$ is that
 of Figure \ref{DC1}. Yet, one would obtain another $\T$-type by changing the lift of the horizontal or of the vertical edge (based at the same top left corner of the square).
 The symmetry of the square makes it sufficient to restrict to the horizontal edge case. Thus, another $\T$-type is represented in
  Figure \ref{ineqDC1}: it
 is the large cell which is filled in by one cell $DC_1$ (of Figure \ref{DC1}), two
cells of type $AA=Br$, and one cell $ABr=BrA$ of level 6 (the small square of Figure \ref{ineqABr=BrA}). The point is that 
all of them belong to $\CR$, so that the boundary of the large cell is homotopically trivial in
 $\CR$.
 
 \begin{figure}
\begin{center}
\includegraphics{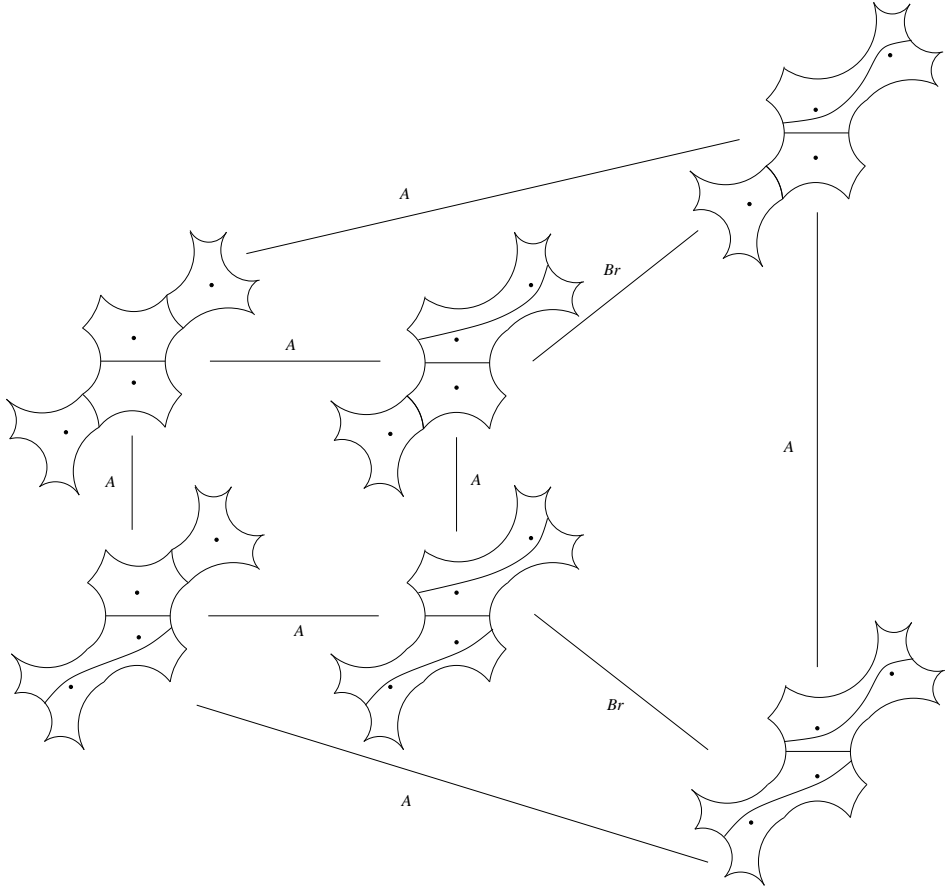}
\caption{Relation between two inequivalent cycles of type $DC_1$}\label{ineqDC1}
\end{center}
 \end{figure}

 \begin{figure}
\begin{center}
\includegraphics{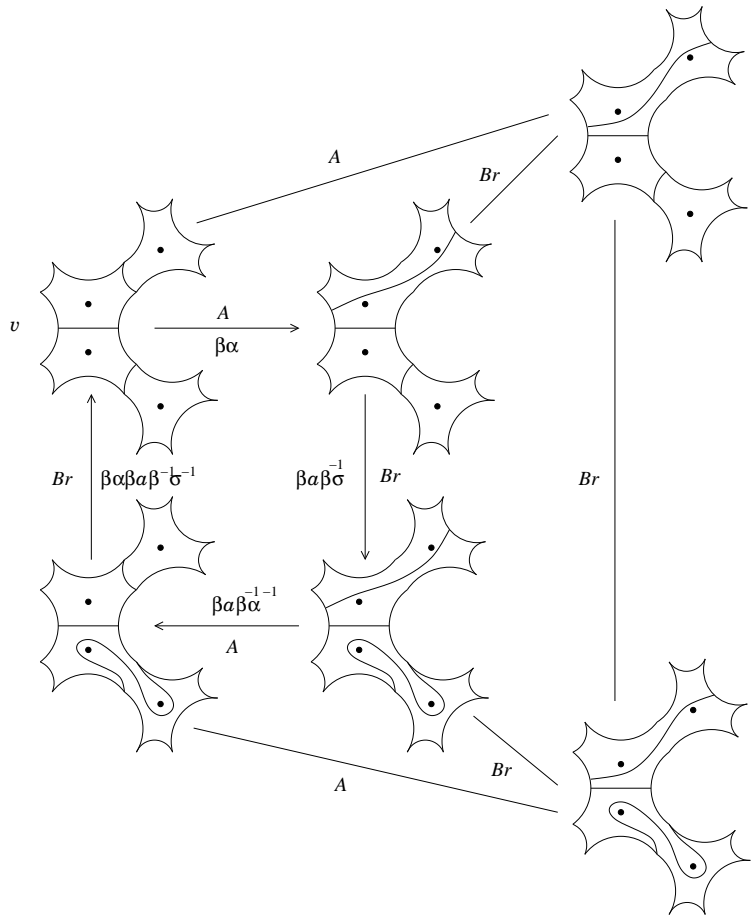}
\caption{Relation between two inequivalent cycles  $ABr=BrA$ of level 6}\label{ineqABr=BrA}
\end{center}
 \end{figure}

\vspace{0.2cm}\noindent 
 Note that another $\T$-type would be obtained by changing the vertical edge of the large cell
 of Figure \ref{ineqDC1}. But we would prove, by the same puzzle game as above, that its
  boundary is trivial in $\CR$. However, there is a subtlety here: the piece $ABr=BrA$ of level 6 we would use would not be $\T$-equivalent to that used above (i.e. the small square of Figure \ref{ineqABr=BrA}). This inequivalent piece $ABr=BrA$ is the large square that is
  represented in Figure \ref{ineqABr=BrA}. But the same figure shows that the latter is equivalent to the cell $ABr=BrA$ (belonging to $\CR$), modulo some cells which all belong
  to $\CR$.  Indeed, the large square
  is filled in by the small piece $ABr=BrA$, two cells $AA=Br$, and one cell $Br_{1}Br_{2}=Br_{2}Br_{1}$
  of level 6. All of them do belong to $\CR$. 
  
  \item  Let us now consider the $\T$-types of cells of $T$-type $DC_{2}$. The only $\T$-type which belongs to $\CR$ is that
 of Figure \ref{DC2}. One would obtain different $\T$-types by changing the lift of the horizontal edge (see the large cell of Figure \ref{ineqDC2}) or of the vertical edge
(see the large cell of Figure \ref{ineqDC2bis}). Using cells of type $DC_{2}$ of Figure \ref{DC2}, of type $AA=Br$, and of type $ABr=BrA$ of level 7
(of Figure \ref{BrA=ABr7} and Figure \ref{BrA=ABr7bis}), as puzzle
 pieces (which all belong to $\CR$), one proves that the large cells of Figures \ref{ineqDC2} and \ref{ineqDC2bis} are equivalent to
 the small ones modulo cells which are in $\CR$.
 
 \begin{figure}
\begin{center}
\includegraphics{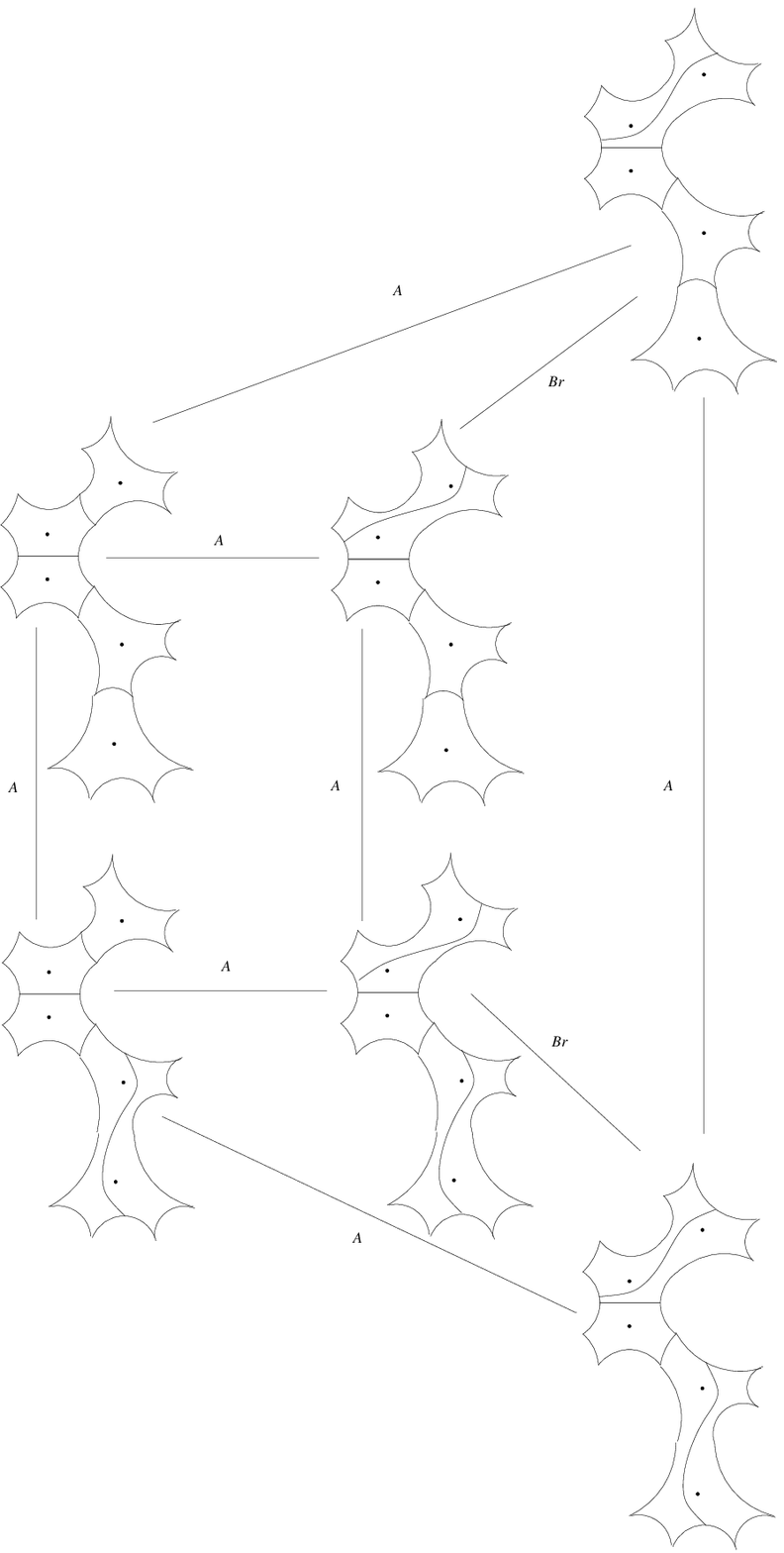}
\caption{Relation between two inequivalent cycles of type $DC_2$}\label{ineqDC2}
\end{center}
\end{figure}

\begin{figure}
\begin{center}
\includegraphics{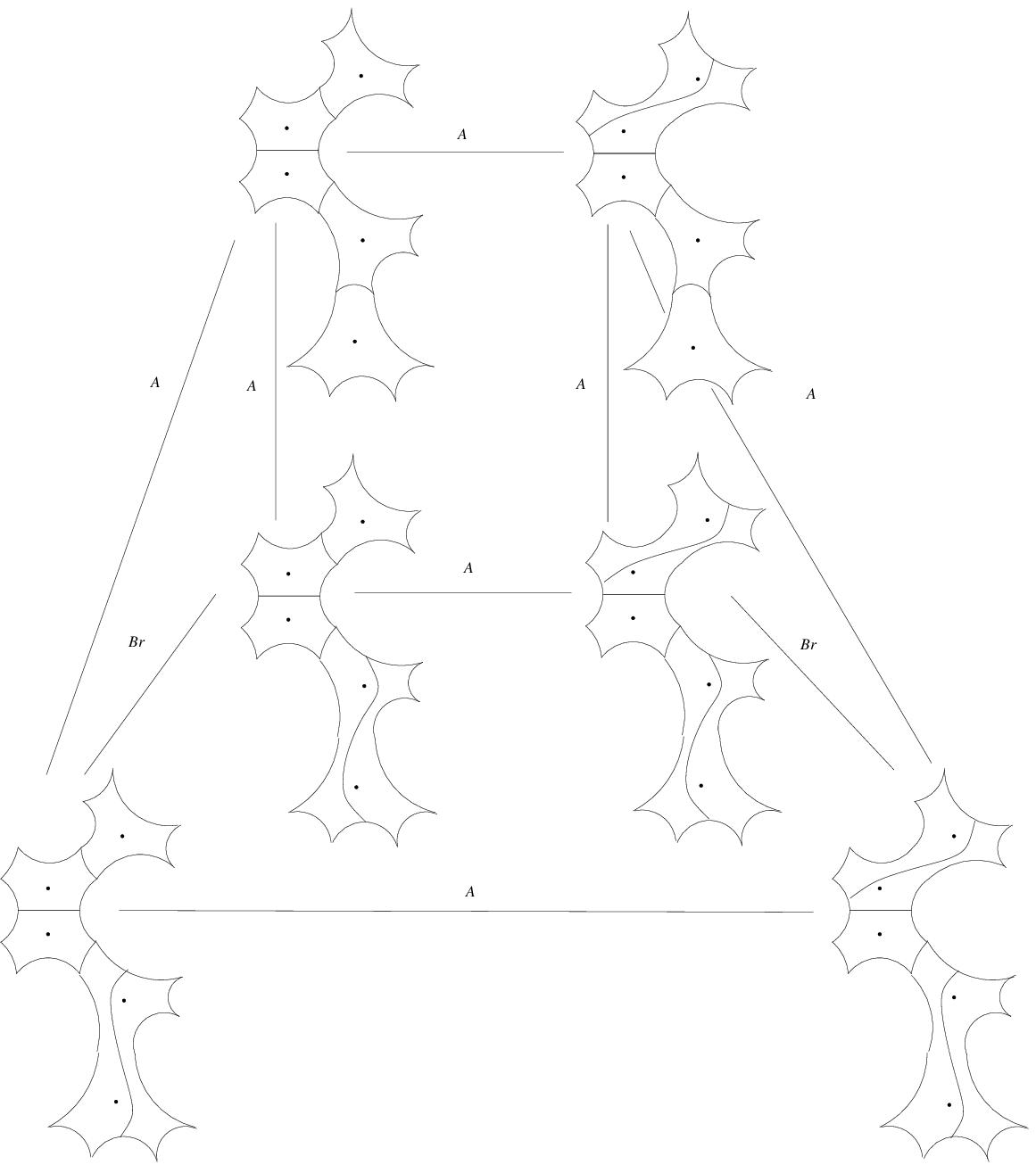}
\caption{Relation between two inequivalent cycles of type $DC_2$, bis}\label{ineqDC2bis}
\end{center}
\end{figure}

\vspace{0.2cm}\noindent 
 Note that one could obtain other $\T$-types of cells of $T$-type $DC_{2}$ by changing both the horizontal
 and the vertical edges. The pieces $ABr=BrA$ of level 7 we would need might not be $\T$-equivalent to those
 of  Figure \ref{BrA=ABr7} and Figure \ref{BrA=ABr7bis}, but equivalent to the latter modulo cells $AA=Br$ and $Br_{1}Br_{2}=Br_{2}Br_{1}$
 of level 7 (Figures \ref{Brlev7} and \ref{Brlev7bis}).
 
 \begin{figure}
\begin{center}
\includegraphics{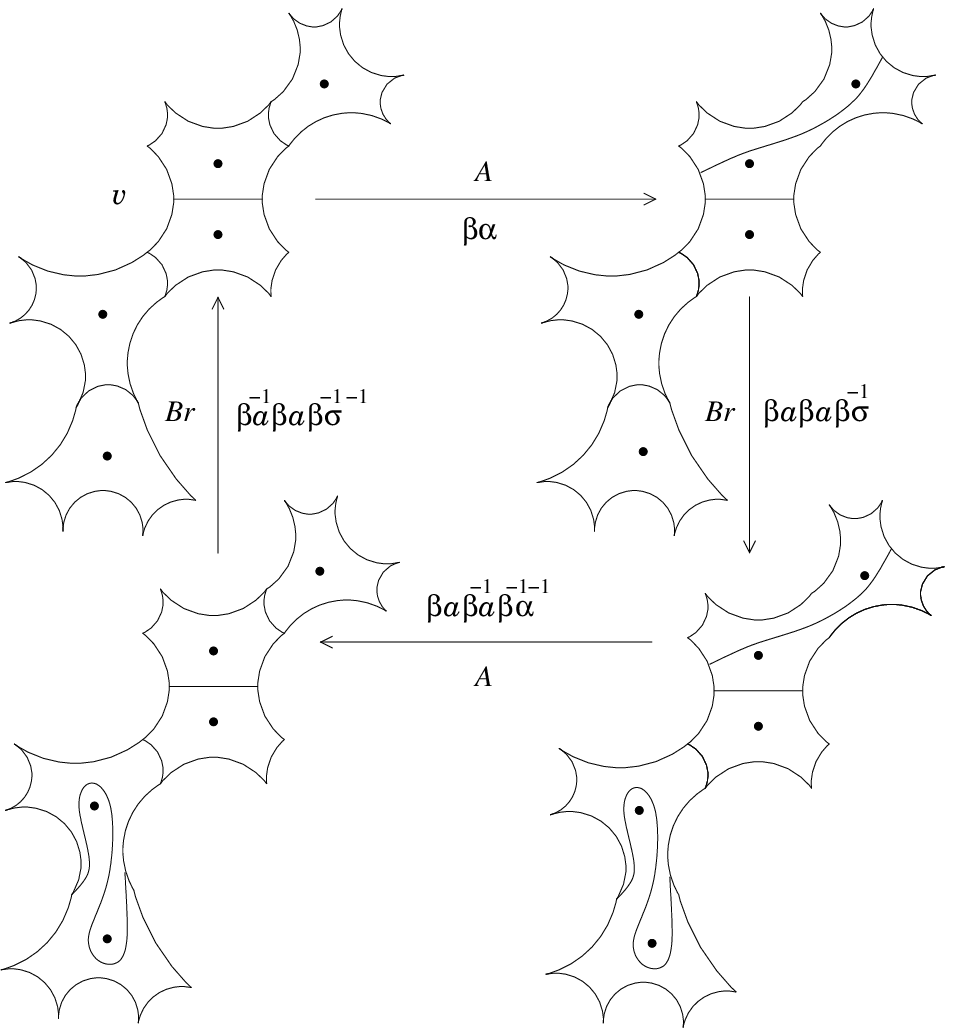}
\caption{First cycle $BrA=ABr$ of level 7}\label{BrA=ABr7}
\end{center}
\end{figure}

\begin{figure}
\begin{center}
\includegraphics{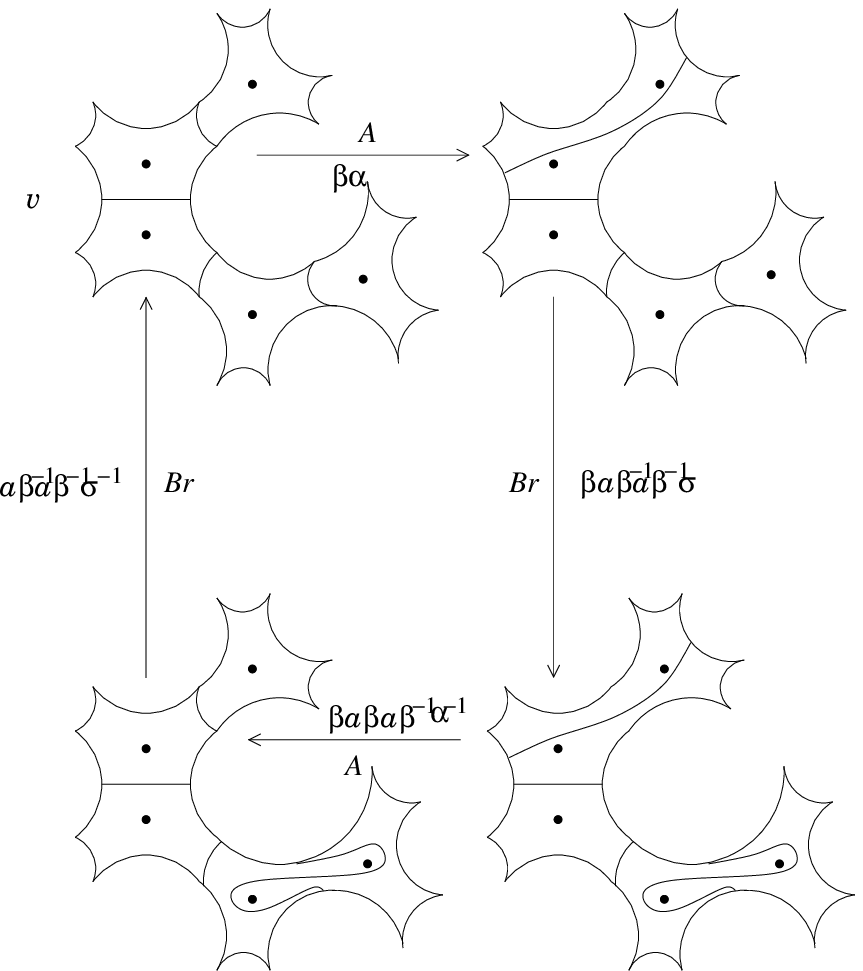}
\caption{Second cycle $BrA=ABr$ of level 7}\label{BrA=ABr7bis}
\end{center}
\end{figure}
 
  \item  Let us finally consider the $\T$-types of cells of pentagonal $T$-type. The only $\T$-type which belongs to $\CR$ is that
 of Figure \ref{penta}. The others (see the large pentagonal cell of Figure \ref{ineqpenta}) are
 equivalent to it modulo cells of type $AA=Br$ and of type $ABr=BrA$ of level 5 (of Figure \ref{BrA=ABr5}), which all belong to $\CR$.
 
 \begin{figure}
\begin{center}
\includegraphics{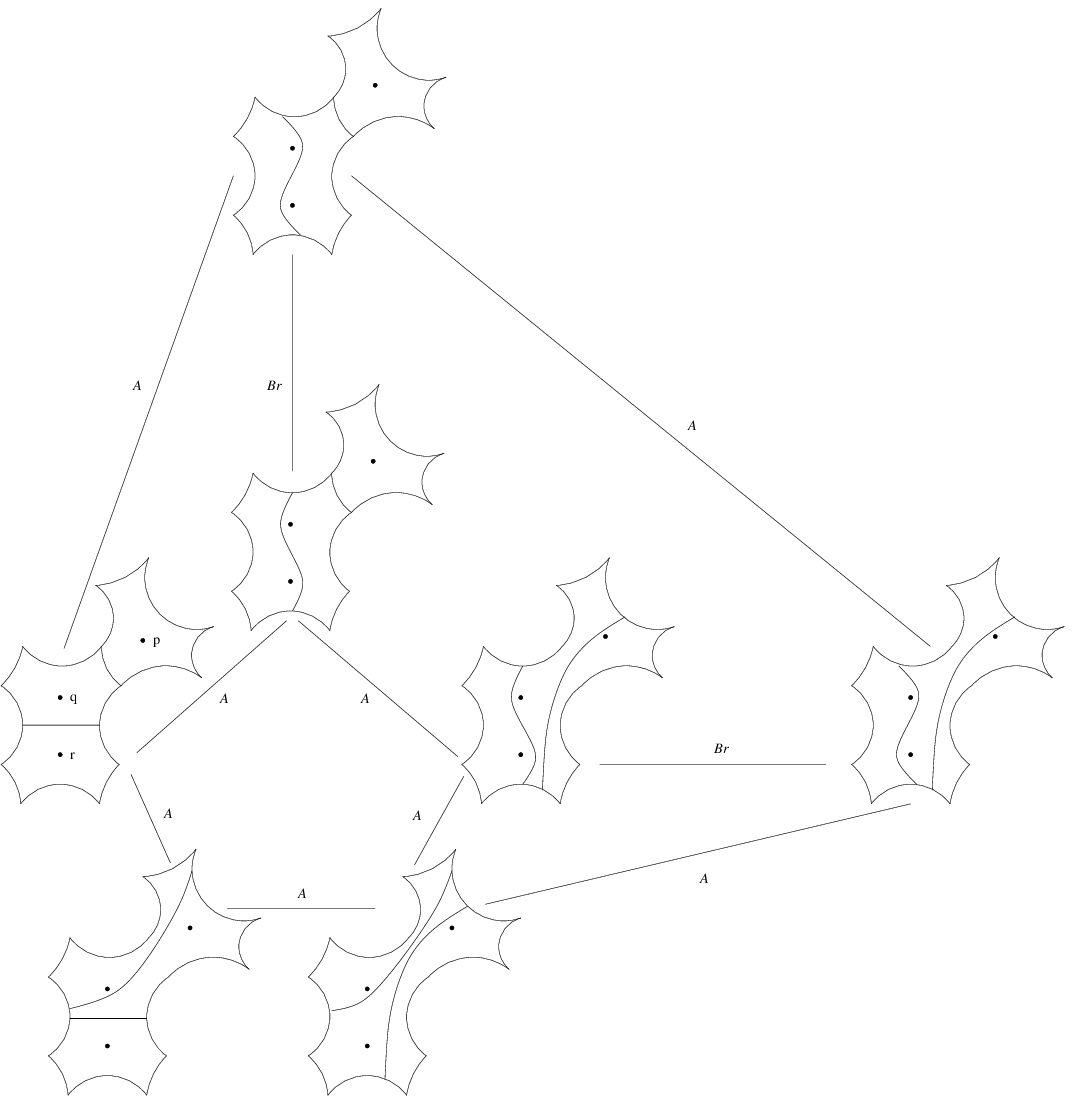}
\caption{Relation between two inequivalent pentagonal cycles}\label{ineqpenta}
\end{center}
\end{figure}
\end{itemize}
\end{proof}
  \vspace{0.2cm}\noindent  
 The second lemma to prove is:

 \begin{lemma}\label{le2}
 Each square cycle of the form ``$A_{1}A_{2}=A_{2}A_{1}$'' in the 1-skeleton of $\C$ or $\CR$, resulting from by the commutation 
 of two $A$-moves along disjoint arcs
  of a rigid structure, may be filled in by 2-cells of $\C$, of $T$-type $DC_{1}$, $DC_{2}$, or of
  pentagonal $T$-type. Therefore, by Lemma \ref{le1}, it may be filled in by 2-cells which all belong to $\CR$, hence is homotopically trivial in
  $\CR$.
 \end{lemma}
 
 \begin{proof}
 Let $\gamma$ denote the square cycle $A_{1}A_{2}=A_{2}A_{1}$ in $\C$ or $\CR$. In \cite{FK}, it is proved that the square $f(\gamma)$
 in  ${\cal HT}_{red}(D)$ may be filled in by 2-cells ($DC_{1}$, $DC_{2}$ and pentagons). Let us enumerate them by $\omega_{1},\ldots\omega_{n}$
 in such a way that $\omega_{i}$ and $\omega_{i+1}$ (for $i=1,\ldots, n-1$) are adjacent along an edge, as well as $\omega_{n}$ 
 and $\omega_{1}$. Following this enumeration, one may lift each $\omega_{i}$ to a 2-cell $\tilde{\omega_{i}}$ of $\C$, in such a
 way that the $n$ lifts fill in the cycle $\gamma$. 
 \end{proof}
 
\vspace{0.2cm} 
\noindent 
 The third and last lemma is:
 
 \begin{lemma}\label{le3}
 The square cycles $Br_1Br_2=Br_2Br_1$ and $ABr=BrA$, bounding the 2-cells of
  $\C$ which are {\em not} in $\CR$, are filled in by some  of 2-cells belonging to $\CR$.  Hence they are homotopically trivial in
  $\CR$.
  \end{lemma}
  
\begin{figure}
\begin{center}
\includegraphics{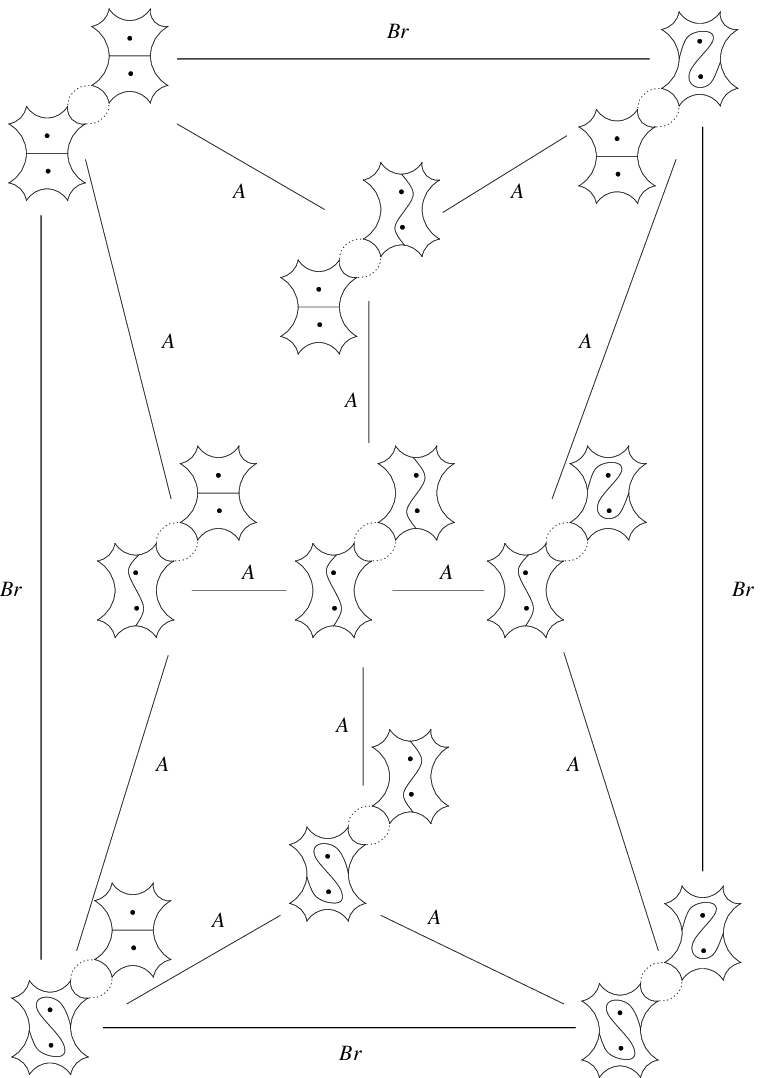}
\caption{A cycle $Br_1Br_2=Br_2Br_1$ is filled in by cells of types $A_1A_2=A_2A_1$ and $AA=Br$}\label{lem1}
\end{center}
\end{figure}

\begin{figure}
\begin{center}
\includegraphics{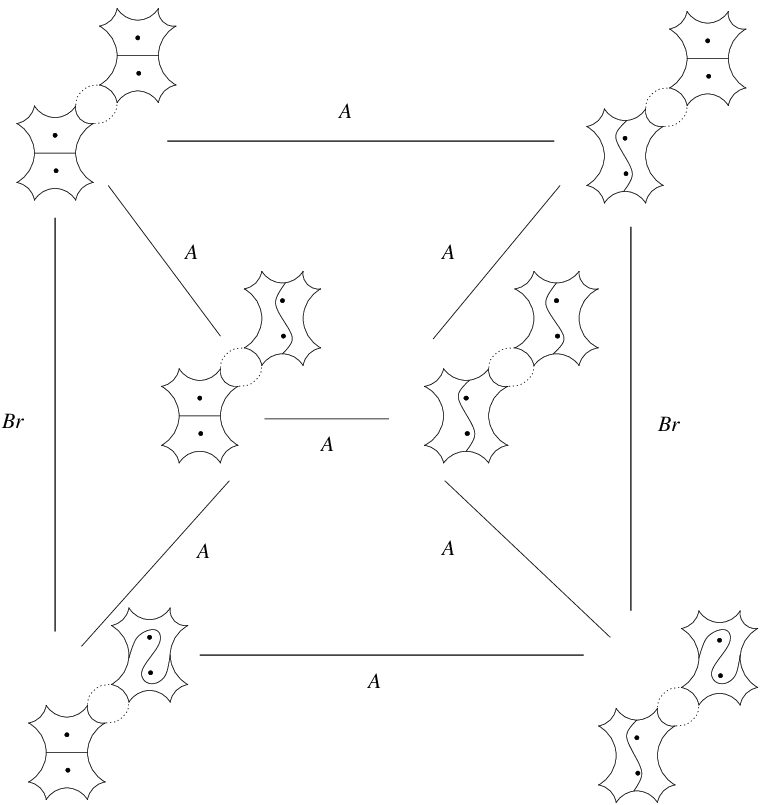}
\caption{A cycle $ABr=BrA$ is filled in by cells of types $A_1A_2=A_2A_1$ and $AA=Br$}\label{lem2}

\end{center}
\end{figure}

  \begin{proof}
 The key point is that a $Br$-move may be seen as the composite of two $A$-moves (``$AA=Br$''), 
 so that each relation of commutation involving $Br$-moves reduces to relations involving $A$-moves. Figure \ref{lem1}
  shows how the square cycles $Br_1Br_2=Br_2Br_1$ are filled in by 4 squares $A_1A_2=A_2A_1$ and four triangles $AA=Br$.
   Since each square $A_1A_2=A_2A_1$  is filled in by cells belonging to $\CR$ by Lemma \ref{le2}, this proves our claim.
    Similarly, Figure \ref{lem2}
    shows how the square cycles $ABr=BrA$ are filled in by 2 cells of type $AA=Br$ and 2 squares $A_1A_2=A_2A_1$.
\end{proof}

\vspace{0.1cm} 
\noindent 
Since the complement of $\CR$ in $\C$ is a union of cells $Br_1Br_2=Br_2Br_1$ and $ABr=BrA$, the last lemma implies that 
the inclusion $\CR\subset\C$ induces an isomorphism at the $\pi_1$ level. This ends the proof of Proposition \ref{main}.   
\end{proof}

\begin{theorem}\label{fp}
The group $\T$ is finitely presented.
\end{theorem}

\begin{figure}
\begin{center}
\includegraphics{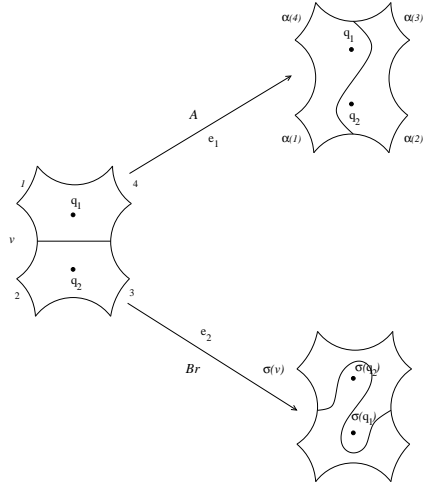}
\caption{Edges $e_{1}$ and $e_{2}$ and generators $\als$ and $\sigma$}\label{rep}
\end{center}
\end{figure}

\begin{proof}
The group $\T$ acts cocompactly on the simply connected complex $\CR$. The stabilizers of the vertices are all isomorphic to $PSL(2,\Z)$. 
Indeed, since all vertices are equivalent modulo $\T$, it suffices to consider the case of the canonical rigid structure $v={\mathfrak r}^*$.
 Clearly, the stabilizer of 
${\mathfrak r}^*$ is also
 the group of orientation-preserving automorphisms of the tree of $D^*$. It is isomorphic to $PSL(2,\Z)$ (cf. also Remark \ref{PSL}).

\vspace{0.2cm} 
\noindent 
We claim that the stabilizers of the edges are isomorphic to $\Z/2\Z$. Indeed, there are two distinct classes of edges (modulo $\T$), which may be represented 
by two edges $e_{1}$ and $e_{2}$, based at the same origin $v$.  The edge $e_1$ corresponds to an $A$-move on an
 arc of reference $\gamma$, while $e_2$ corresponds to a $Br$-move on the same arc $\gamma$. We may assume that the mapping
  class $\als\in\T$ of \ref{2el} has been chosen such that the terminal vertex of $e_1$ is $\als(v)$. Recall
 that it is
  a rigid rotation of order 4 outside $H_1\cup H_2$, but it fixes $q_1$ and $q_2$ inside $H_1\cup H_2$, see Figure \ref{rep}.

\vspace{0.2cm} 
\noindent 
  We denote by  $\sigma\in B_{\infty}\subset \T$ the positive braiding on the arc $\gamma$.  It permutes 
  $q_1$ and $q_2$ and is such that $\sigma(v)$ is the terminal  vertex of $e_2$, see Figure \ref{rep}.

\vspace{0.2cm} 
\noindent 
   One first checks that
 there is no element of $\T$ that reverses the orientation of the edges $e_{1}$ or $e_{2}$. Thus, the stabilizers $\T_{e_{1}}$ of $e_{1}$ and
 $\T_{e_{2}}$ of $e_{2}$
 are subgroups of the stabilizer of $v$. In fact, $\T_{e_{1}}=\T_{e_{2}}$, generated by the element $a\in PSL(2,\Z)$ of order 2, 
 which is rigid rotation of angle $\pi$ that interchanges the hexagons $H_{1}$ and $H_{2}$ and  preserves the arc $\gamma$ (reversing its orientation).
 We shall see below that $a= \sigma^{-1}\als^2= \als^2\sigma^{-1}$  (beware that ${\alpha^*}^2$ is not of order 2!).

\vspace{0.2cm} 
\noindent 
 Since the stabilizers of the vertices are finitely presented and the stabilizers of the edges are finitely generated, Theorem 1 of \cite{br}
 asserts that $\T$ is finitely presented.
\end{proof} 

\begin{remark}
The stabilizer  $T^*_v$ of $v$ admits the following presentation:
$$\T_v=<a,\bs\,|\, a^2=\bs^3=1>\cong PSL(2,\Z).$$
\end{remark}

\subsection{A presentation for $\T$}

\subsubsection{Statement of the theorem}

\begin{theorem}

The group $\T$ admits a finite presentation, with three generators $\als$, $\bs$ and $\sigma$, and the following relations:

\begin{enumerate}
\item $\als\sigma=\sigma\als$
\item $\als^4=\sigma^2$
\item $\bs^3=1$
\item $(\bs\als)^5=\sigma \bs\sigma\bs^{-1}\sigma$\\

Setting $a=\als^2\sigma^{-1}$,
\item $[\bs\als\bs,a\bs\als\bs a]=1$ (level 5)
\item $[\bs\als\bs, a\bs a\bs\als\bs a\bs a]=1$\\

Consistency relations: 
\item $[\sigma,\bs\als\bs]=1$ (level 5)
\item $[\sigma,\bs\als\bs a\bs^{-1}]=1$ (level 6)
\item $[\sigma,\bs a\bs^{-1}a\bs\als\bs a\bs a\bs^{-1}]=1$ (level 7)
\item $[\sigma,\bs a\bs a\bs\als\bs a\bs a\bs]=1$ (level 7) \\

Commutations of braidings:
\item $[\sigma,\bs a\bs^{-1}\sigma (\bs a\bs^{-1})^{-1}]=1$ (level 6)
\item $[\sigma,\bs a\bs a\bs^{-1}\sigma (\bs a\bs a\bs^{-1})^{-1}]=1$ (level 7)
\item $[\sigma,\bs a\bs a\bs\sigma (\bs a\bs a\bs )^{-1}]=1$ (level 7)\\

Setting $\sigma_1=\sigma$, $\sigma_2=\bs\sigma\bs^{-1}$ and $\sigma_3=\bs^{-1}\sigma\bs$,
\item $\sigma_1\sigma_2\sigma_1= \sigma_2\sigma_1\sigma_2$ (fundamental relation of the braid group)
\item $\sigma_1\sigma_2\sigma_3\sigma_1=\sigma_2\sigma_3\sigma_1 \sigma_2=\sigma_3\sigma_1 \sigma_2 \sigma_3$ (Sergiescu's relations)
\end{enumerate}

\end{theorem}

\begin{corollary}
We have $H_{1}(\T)=\Z/6\Z$. In particular, the groups $\T$ and $T^{\star}$ are not isomorphic.
\end{corollary}

\begin{proof}
 $H_{1}(\T)$ is generated by the commuting  $[\als]$,  $[\bs]$ and $[\sigma]$, subject to the relations $4[\als]=2[\sigma]$, $3[\bs]=0$ and $5[\als]+
 5[\bs]=3[\sigma]$. They are equivalent to $[\als]=-[\sigma]$, $[\bs]=-2[\sigma]$ and $6[\sigma]=0$, hence the claim.
 \end{proof}

\subsubsection{Generators}

We follow the method described by K. Brown in \cite{br}, derived from the Bass-Serre theory.

\vspace{0.2cm} 
\noindent 
Recall that the quotient $\CR/\T$ possesses a unique vertex, represented by $v$, and two edges, represented by $e_{1}$
 and $e_{2}$ (see the proof of
Theorem \ref{sta}).

\vspace{0.2cm} \noindent 
The stabilizers of $e_1$ and $e_2$ are of order 2, generated by $a$. For better clarity, we write $\T_{e_1}=<a_1,\; a_{1}^2=1>$ and $\T_{e_2}=<a_2,\; a_{2}^2=1>$, respectively.

\subsubsection{Relations}

Theorem 1 of \cite{br} states that  $\T$ is generated by the  
stabilizer $\T_v$ and by the elements $g_{e_1}=\als$ and $g_{e_2}=\sigma$, subject to the         
following relations:         
\begin{itemize}         
\item {\it Pres. (i):} For each $e\in \{e_1,e_2\}$, $g_e^{-1} i_e(h) g_e=c_e(h)$ for    
  all $h\in \T_e$,          
where $i_e$ is the                       
inclusion $\T_e\hookrightarrow \T_{v}$ and $c_e:\T_e\rightarrow \T_v$          
is the conjugation morphism $h\mapsto  g_e^{-1} h  g_e$.

\vspace{0.1cm} 
\noindent 
Explicitly, ``$g_{e_1}^{-1} i_{e_1}(h) g_{e_1}=c_{e_1}(h)$'' with  $g_{e_1}=\als$ and $h=a_1\in \T_{e_1}$ provides the relation  $\als^{-1}a\als=\als^{-1}a_1\als$, where the right hand side is computed in $\T_v$, in which it is equal to $a$. Hence the relation
$$\als^{-1}a\als=a \;\;(*).$$

\vspace{0.1cm} 
\noindent 
As for ``$g_{e_2}^{-1} i_{e_2}(h) g_{e_2}=c_{e_2}(h)$'' with  $g_{e_2}=\sigma$ and $h=a_2\in \T_{e_2}$, it provides the relation $\sigma^{-1}a\sigma= \sigma^{-1}a_2\sigma$, where the right hand side is computed in $\T_v$, in which it is equal to $a$. Hence the relation

$$\sigma^{-1}a\sigma=a\;\;(**).$$

\item {\it Pres. (ii):}  $r_{\tau}=1$ for each 
$2$-cell $\tau\in \CR$,          
where $r_{\tau}$ is a word in the generators of $\T_v$, $\als$ and          
$\sigma$, associated with the $2$-cell $\tau$ in the way described in \cite{br}. We recall it for the convenience of the reader: 
         
\end{itemize}

\vspace{0.1cm} 
\noindent  
Each edge of the complex starting at $v$ has one of the following forms: 
 
\begin{enumerate} 
\item $v\longrightarrow (h\als^{\pm 1})(v)$, $h\in\T_v$
\item $v\longrightarrow (h\sigma^{\pm 1})(v)$, $h\in\T_v$
\end{enumerate}

\vspace{0.1cm} 
\noindent 
To such an edge $e$ we associate an element $\gamma\in \T$ such that $e$ 
ends at $\gamma(v)$: $\gamma=h\als^{\pm 1}$ in case (a), $\gamma=h\sigma^{\pm 1}$ in case (b).

\vspace{0.1cm} 
\noindent 
Let $\tau$ be one of the  2-cells of the complex $\CR$. One chooses an orientation and a cyclic 
labeling of the boundary edges, such that the labeled 1 edge $E_1$  starts from 
the vertex $v$.

\vspace{0.1cm} 
\noindent
Let $\gamma_1$ 
be associated to $E_1$ as above. It ends at $\gamma_1(v)$, so the second edge 
is of the form $\gamma_1(E_2)$ for some edge $E_2$ starting at $v$. Let 
$\gamma_2$ be associated to $E_2$. The second edge ends at 
$\gamma_1\gamma_2(v)$. If $n$ is the length of the cycle bounding $\tau$, one obtains this way a sequence $\gamma_1,\ldots,\gamma_n$ such that 
$\gamma_1\cdots\gamma_n(v)=v$.

\vspace{0.1cm} 
\noindent
Note that for each of the  cycles, we have 
indicated the corresponding $\gamma_i$ above the $i^{th}$ edge.

 \vspace{0.1cm} 
\noindent 
Let $\gamma$ be the element of the stabilizer $\T_v$ which is equal to $\gamma_1\cdots\gamma_n$ when each element $\gamma_i$ is viewed in $\T$. Then 
the relation associated to $\tau$ is 
$$ \gamma_1\cdots\gamma_n=\gamma$$ 
where the left hand side is viewed as a word in $\als$, $\sigma$, $a$, $\bs$, and their inverses.

\vspace{0.1cm} 
\noindent 
Following this process for the 2-cells of the complex $\CR$, one obtains:

\vspace{0.2cm} 
\noindent {\it 1. Cell $AA=Br$ (Figure \ref{AA=Br})}. The corresponding relation is $\als\sigma^{-1}\als=a$. Equivalently, 
$\sigma=\als a^{-1} \als$. Since by $(*)$, $\als$ and $a$ commute, one obtains $\sigma=\als^2 a^{-1}$, 
hence $\als$ and $\sigma$ commute ({\it Rel. 1}. Since $a=\sigma^{-1}\als^2$, $a$ may be eliminated,
 and the relation $a^2=1$ is now equivalent to $\als^4=\sigma^2$ ({\it Rel. 2}).

\vspace{0.2cm}  
\noindent {\it 2. a) Cell  $DC_{1}$ (Figure \ref{DC1})}. The corresponding relation is $[\bs\als\bs,a\bs\als\bs a]=1$ 
({\it Rel. 5}).

\vspace{0.2cm} 
\noindent {\it b) Cell  $DC_{2}$ (Figure \ref{DC2})}. The corresponding relation is  $[\bs\als\bs, a\bs a\bs\als\bs a\bs a]=1$ ({\it Rel. 6}).

\vspace{0.2cm} 
\noindent {\it c) Pentagonal cell} (Figure \ref{penta}).
It gives first the relation $\als^{-1}\bs^{-1}\als^{-1}\bs^{-1}\als^{-1}\sigma\bs^{-1}\als^{-1}\sigma\bs^{-1}\als=\bs a$. Taking
 the inverse of this relation, one obtains
 $\als^{-1}\bs\sigma^{-1}\als\bs\sigma^{-1}\als\bs\als\bs\als=a\bs^{-1}$.
 Replacing $a$ by $a=\als^{-2}\sigma$, one obtains $\als\bs\sigma^{-1}\als\bs\sigma^{-1}(\als\bs)^{3}=\sigma$.
 Equivalently, $\sigma^{-1}\als\bs\sigma^{-1}(\als\bs)^{4}=(\als\bs)^{-1}\sigma\als\bs$. Since $\als$
 and $\sigma$ commute, the right hand side is equal to $\bs^{-1}\sigma\bs$. Hence the relation 
 $\sigma^{-1}\als\bs\sigma^{-1}(\als\bs)^{-1}(\als\bs)^{5}=\bs^{-1}\sigma\bs$. This is equivalent to
 $(\bs\als)^{5}=\bs\sigma\bs^{-1}\als^{-1}\sigma\bs^{-1}\sigma\bs\als$. Using once again the commutation between $\sigma$
 and $\als$, one obtains 
 $(\bs\als)^{5}=\bs\sigma\bs^{-1}\sigma\als^{-1}\bs^{-1}\sigma\bs\als$. But we shall see below that $\sigma$ and
 $\bs\als\bs$ commute (cf. {\it Cell $ABr=BrA$ of level 5}), so that $\als^{-1}\bs^{-1}\sigma\bs\als=\bs\sigma\bs^{-1}$. Finally,
 the relation becomes
 $$(\bs\als)^{5}=\bs\sigma\bs^{-1}\sigma \bs\sigma\bs^{-1}.$$
 Modulo the braid relation (cf. below {\it Cells coming from the presentation of the braid group}), this is {\it Rel. 4}.

\vspace{0.2cm} 
\noindent {\it 3. a) Cell $ABr=BrA$ of level 5 (Figure \ref{BrA=ABr5})}:  $[\sigma,\bs\als\bs]=1$ ({\it Rel. 7}).

\vspace{0.2cm} 
\noindent {\it b) Cell $ABr=BrA$ of level 6 (Figure \ref{ineqABr=BrA})}: $[\sigma,\bs\als\bs a\bs^{-1}]=1$ ({\it Rel. 8}).

\vspace{0.2cm} 
\noindent {\it c) First cell $ABr=BrA$ of level 7 (Figure \ref{BrA=ABr7})}: $[\sigma,\bs a\bs^{-1}a\bs\als\bs a\bs a\bs^{-1}]=1$ ({\it Rel. 9}).

\vspace{0.2cm} 
\noindent {\it d) Second cell $ABr=BrA$ of level 7 (Figure \ref{BrA=ABr7bis})}:  $[\sigma,\bs a\bs a\bs\als\bs a\bs a\bs]=1$  ({\it Rel. 10}).

\vspace{0.2cm} 
\noindent {\it 4. a) Cell $Br_{1}Br_{2}=Br_{2}Br_{1}$ of level 6 (Figure \ref{Brlev6})}: $[\sigma,\bs a\bs^{-1}\sigma (\bs a\bs^{-1})^{-1}]=1$ ({\it Rel. 11
}).

\vspace{0.2cm} 
\noindent {\it b) First cell $Br_{1}Br_{2}=Br_{2}Br_{1}$ of level 7 (Figure \ref{Brlev7})}: $[\sigma,\bs a\bs a\bs^{-1}\sigma (\bs a\bs a\bs^{-1})^{-1}]=1$
 ({\it Rel. 12}).

\vspace{0.2cm} 
 \noindent {\it c) Second cell $Br_{1}Br_{2}=Br_{2}Br_{1}$ of level 7 (Figure \ref{Brlev7bis})}: $[\sigma,\bs a\bs a\bs\sigma (\bs a\bs a\bs )^{-1}]=1$
({\it Rel. 13}).

\vspace{0.2cm} 
\noindent {\it 5. Cells coming from the presentation of the braid group}. They obviously give the relations {\it 14} and {\it 15}.

\subsection{$T^{\star}$ is finitely presented}

The groups $\T$ and $T^{\star}$, though both alike, are not isomorphic. However, there is a proof for the assertion that $T^{\star}$ 
is finitely presented which mimics that for $\T$. One introduces $T^{\star}$-complexes ${\cal C }^+(T^{\star})$ and 
${\cal C }(T^{\star})$, whose vertices are the asymptotically rigid structures of $D^{\star}$, and the edges are of two types, corresponding to moves 
$A$ and $Br$.

\begin{figure}
\begin{center}
\includegraphics{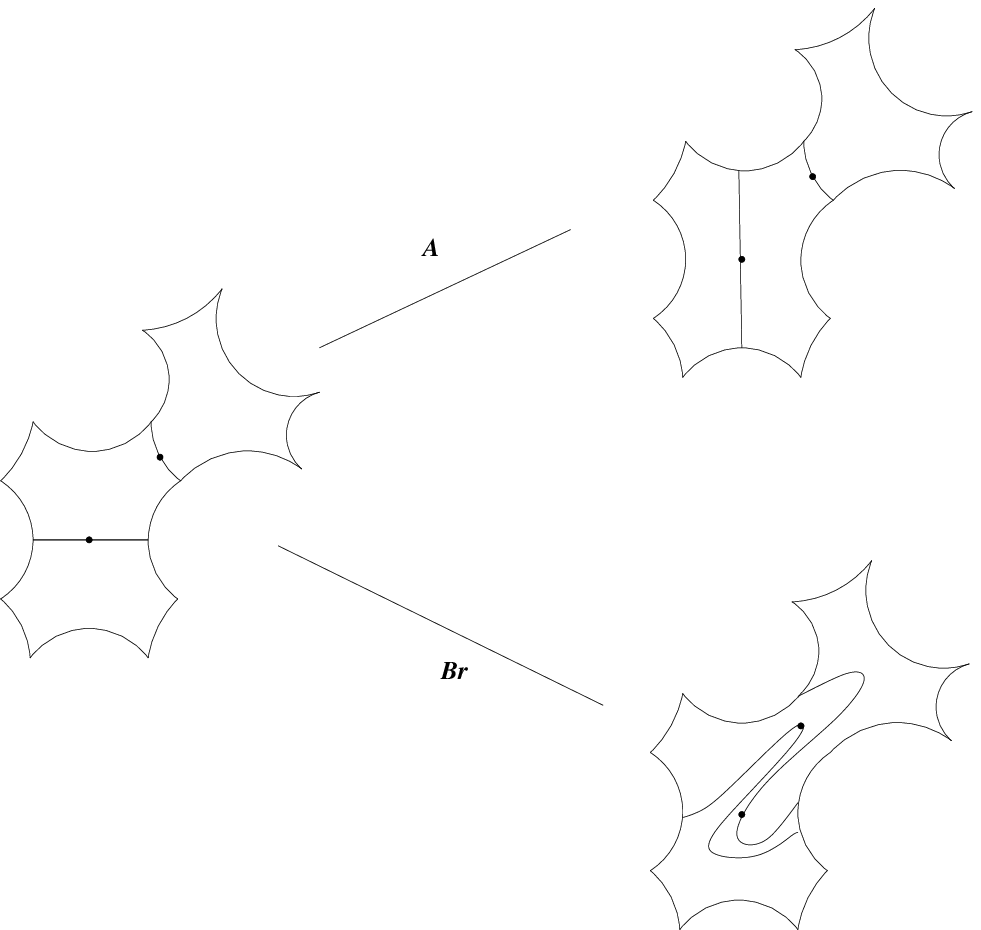}
\caption{The two types of edges in ${\cal C}(T^{\star})$}\label{mouv2}
\end{center}
\end{figure}

\begin{itemize}
\item  If  $\mathfrak{r}$ is an asymptotically rigid structure and $\gamma$  is an arc of $\mathfrak{r}$, the $A$-move on  $\gamma$
keeps unchanged all the arcs of $\mathfrak{r}$ except $\gamma$, and replaces $\gamma$ by $\gamma'$ which is transverse to $\gamma$ and 
passes through the same puncture
as $\gamma$ (see Figure \ref{mouv2}). Note that there is a unique $A$-move on $\gamma$.

\item 
If $\mathfrak{r}$ is an asymptotically rigid structure and $p$ and $q$ are two punctures of $D^{\star}$ on two 
different sides of a hexagon $H$ of $\mathfrak{r}$, there is a simple arc $e$ inside $H$ which connects $p$ to $q$. Let $\sigma_{e}$ be
the braiding along $e$. The move $Br$ changes $\mathfrak{r}$ by the natural action of $\sigma_{e}$ on  $\mathfrak{r}$
  (see Figure \ref{mouv2}).
\end{itemize}

\begin{figure}
\begin{center}
\includegraphics{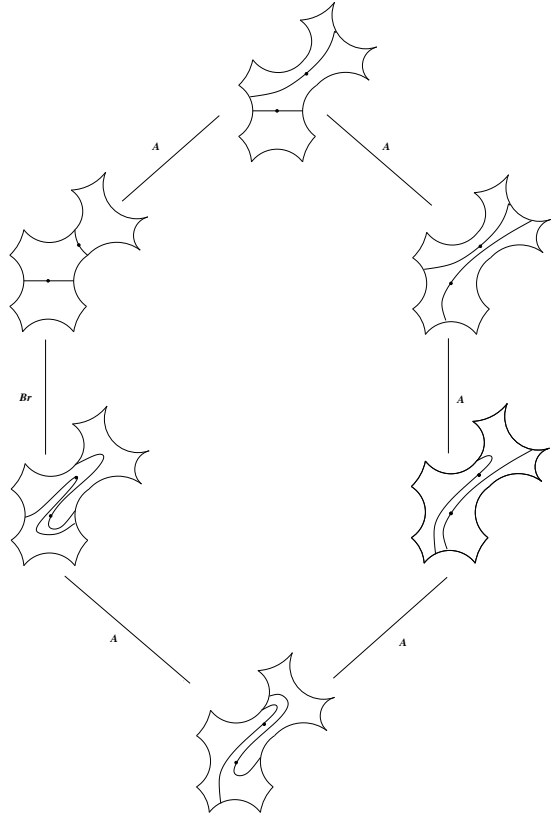}
\caption{Relation $``A^5=Br"$}\label{pentaBr}
\end{center}
\end{figure}

\vspace{0.2cm} \noindent 
The complex  ${\cal C }^+(T^{\star})$ has the same types of 2-cells as ${\cal C }^+(\T)$, except that:\\
\noindent -- the cell $AA=Br$ 
of ${\cal C }^+(\T)$ does not exist, \\
-- and the pentagonal cell is replaced by a hexagonal cell, expressing that a certain sequence of five $A$-moves produces the
effect of a $Br$-move, see Figure \ref{pentaBr}. As can be guessed, the relation which will be associated to this cell is $(\beta^{\star}\alpha^{\star})^{5}=
\sigma_{[04]}$, 
see Lemma \ref{sigmapenta}.

\vspace{0.2cm}\noindent  
The Cayley complex of $B_{\infty}$, for the Sergiescu presentation associated to the graph of $D^{\star}$ (cf. Remark \ref{graph}), 
provides 2-cells of ${\cal C}(T^{\star})$ which involve only $Br$-moves. There are infinitely many $T^{\star}$-types 
of cells $Br_{1}Br_{2}=Br_{2}Br_{1}$, but since the graph is homogeneous, the other cells of the Cayley complex of $B_{\infty}$ provide
finitely many other cells in ${\cal C}(T^{\star})/T^{\star}$. Therefore, just like ${\cal C}(\T)/\T$, the quotient ${\cal C}(T^{\star})/T^{\star}$ is
made of countably many 2-cells $Br_{1}Br_{2}=Br_{2}Br_{1}$ and $ABr=BrA$, plus finitely many other 2-cells.

\vspace{0.2cm} \noindent
Fortunately, the key role played by the 2-cell $AA=Br$ in the proof of Theorem \ref{fp}, especially in lemma \ref{le3}, 
is now played by 
the 2-cell $A^5=Br$: one eliminates almost all the $T^{\star}$-types of 2-cells $ABr=BrA$ and $Br_{1}Br_{2}=Br_{2}Br_{1}$ using cycles
$A_{1}A_{2}=A_{2}A_{1}$ (the analogue of Lemma \ref{le2} is true) and $Br=A^5$. This enables us to obtain a reduced complex ${\cal C }(T^{\star})$ which is finite modulo 
$T^{\star}$ and simply connected. These arguments constitute a sketch of the proof of

\begin{theorem}
The braided Ptolemy-Thompson group $T^{\star}$ is finitely presented, and admits a presentation with 2 generators.
\end{theorem}

\section{Comments and open questions}

\noindent {\bf Actions by homeomorphisms on $S^1$.}
D.Calegari (\cite{Ca})
proved that punctured mapping class groups have a faithful action 
by homeomorphisms on $S^1$. Specifically, 
let $S$ be a surface (possibly of infinite type) with
a base point $p$. Let ${\rm Mod}(S,p)$ be the group
${\rm Homeo}^+(S,p)/{\rm Homeo}^+_0(S,p)$. 
Here ${\rm Homeo}^+(S,p)$ denotes the group of orientation-preserving
homeomorphisms of $S$ to itself which takes $p$ to itself, and
${\rm Homeo}^+_0(S,p)$ denotes the connected subgroup containing the
identity map. Then ${\rm Mod}_{S,p}$ is circularly orderable.
Notice that this punctured mapping class groups fits into an exact 
sequence 
\[ 1 \to \pi_1(S) \to {\rm Mod}(S,p) \to {\rm Mod}(S) \to 1\]
where ${\rm Mod}(S)$ is the usual mapping class group of $S$.

\vspace{0.2cm}
\noindent 
In particular, there are  extensions by free groups  
of the universal mapping class 
groups $\B$ in genus zero, which  embed
in ${\rm Homeo}_+(S^1)$, being circularly orderable. 
It seems, however, that $\B$ does not act  faithfully on the circle. 

\vspace{0.2cm}
\noindent {\bf Smoothing the action.}
It is presently unknown whether the group $T^*$ admits an embedding into 
some group of piecewise linear homeomorphisms $PL_{\lambda}(S^1)$ with 
break points and derivatives  of the form $\lambda^n$,
where $n\in\Z$.
We conjecture that there is no embedding 
$T^*\to {\rm Diff}^2_+(S^1)$ into the  group of diffeomorphisms  
of $S^1$ of class ${\mathcal C}^2$. Specifically, any 
homomorphism $T^*\to {\rm Diff}^2_+(S^1)$ should factor through 
a finite extension of $T$.

\vspace{0.2cm}
\noindent {\bf Automatic groups.}
Thompson groups are known to be asynchronous automatic groups
(see \cite{GN}), but it is still unknown whether they are 
(synchronously) automatic. V.Guba proved that the Dehn function 
of $F$ is quadratic, as is the case for all automatic groups. 
We conjecture that $T^*$ is automatic. 
Notice that braid groups and more generally 
mapping class groups are known to be automatic (see \cite{ECHLPT,Mo}).  
In particular, this would immediately imply that $T^*$ is finitely presented and has solvable word 
problem. We expect that the conjugacy problem is  
solvable too, though it is presently unknown whether this holds true 
for all automatic groups. Moreover, automatic groups are combable 
(see \cite{ECHLPT}, p.84) and hence they are 
${\rm FP}_{\infty}$ and thus ${\rm F}_{\infty}$, i.e. they have a 
classifying space with finitely many cells in each dimension
(see \cite{ECHLPT}, p.220). Eventually, the Dehn function of $T^*$ should be quadratic.  If $T^*$ is biautomatic then it would
provide an example of such a group having a free abelian subgroup 
of infinite rank.

\vspace{0.2cm}
\noindent {\bf Outer automorphisms groups.}
It was first established by Dyer and Grossman that 
${\rm Out}(B_n)=\Z/2\Z$, for $n\geq 4$, and Ivanov and further McCarthy
extended this to mapping class groups. 
However, their result does not extend, as stated, to infinite braid 
groups. In fact there exists an embedding 
\[ T\to {\rm Out}(B_{\infty})\]
induced by the action of $T^*$ by conjugacy on its normal subgroup 
$B_{\infty}$.  In particular, ${\rm Out}(B_{\infty})$ seems to be a
quite rich group. 

\vspace{0.2cm}
\noindent 
On the other hand  M.Brin (\cite{Bri0}) 
proved that group ${\rm Out}(T)=\Z/2\Z$. 
It would be interesting to know whether ${\rm Out}(T^*)=\Z/2\Z$ holds. 

\vspace{0.2cm}
\noindent {\bf Spherical generalization: the group $\B^*$.} Let $\mathscr{S}^*$ be the surface obtained by gluing
together the surfaces  $D$ and $D^*$ along their boundaries. A homeomorphism of  $\mathscr{S}^*$ is
 {\it asymptotically rigid} if it maps almost every hexagon of $D$ onto a hexagon of $D$, and almost every hexagon of $D^*$ 
 onto a hexagon of $D^*$. Here ``almost every'' means all but finitely 
many of them. 
 \begin{definition}
 The group ${\cal B}^*$ is the group of isotopy classes of asymptotically rigid homeomorphisms of $\mathscr{S}^*$.
 \end{definition} 
 
\noindent  Let $V$ denote the Thompson group acting on the Cantor set. 
There is  then a short exact sequence
 $$1\rightarrow PM(\mathscr{S}^*)\longrightarrow {\cal B}^*\longrightarrow V\rightarrow 1,$$
 where $PM(\mathscr{S}^*)$ denotes the  pure mapping class group 
of $\mathscr{S}^*$ i.e. the compactly supported mapping classes 
of homeomorphisms of $\mathscr{S}^*$ that preserve the ends 
of $\mathscr{S}^*$. 
 
\vspace{0.1cm}
\noindent As alluded in the introduction, the methods of section 3 can be adapted to prove that 
\begin{proposition}
${\cal B}^*$ is a finitely 
 presented group.
\end{proposition}

\vspace{0.2cm}
\noindent {\bf Remarks.}
 
\begin{enumerate}
\item Any diagram group (see \cite{GuS}) can be 
embedded into $B_{\infty}$ (by a result of Wiest in \cite{Wi}) and 
thus into $T^*$. However, 
$T^*$ is not a diagram group since it has torsion. Moreover,
$T^*$ and the Brin-Dehornoy braided Thompson group $BV$ 
are the typical examples of some more general braided diagram groups.
The work of Farley on  diagram groups 
extends to braided diagram (see \cite{Fa}). In particular,  
each one of these groups acts properly cellularly on a 
$CAT(0)$-complex, which is not locally finite. The stabilizers of cells 
are isomorphic to braid groups (on finitely many strands).  
\item The group $T^*$ has not the Kazhdan property since $T$ has not. 
Moreover it is a-T-menable, by the same reason. 
\item If $\Gamma$ is a lattice in a simple Lie group of 
and the $\Q$-rank of $\Gamma$ 
is at least 2 then any homomorphism $\Gamma\to T^*$ should be trivial, 
since any ${\mathcal C}^0$-action of such a $\Gamma$ on $S^1$ is trivial, 
by a result of D.Witte (see \cite{Wit}). 
\item There exist however homomorphisms from arithmetic groups 
of rank one into $T^*$. In fact, Kontsevich and Soibelman recently constructed 
in \cite{KoS}  faithful  homomorphisms from an arithmetic 
subgroup of $SO(1,18)$ into the braid groups. 
\item The group $T^*$ is non-amenable and hence of exponential growth.  

\end{enumerate}

{
\small      
      
\bibliographystyle{plain}

\begin{thebibliography}{30}      


\bibitem{Bez}
V.N. Bezverkhnii, 
{\em Solution of the generalized conjugacy problem for words in
 $C(p)\&T(q)$-groups}, 
Izv. Tul. Gos. Univ. Ser. Mat. Mekh. Inform. 4 (1998), 
no. 3, Matematika, 5-13.     

      
\bibitem{BK2}      
B. Bakalov and A.Kirillov Jr.,      
{\em On the Lego-Teichmuller game},       
Transform. Groups 5(2000),  207-244.       
   
 
\bibitem{BKL}
J.Birman, Ki Hyoung Ko and  Sang Jin Lee, {\em  A new approach to the word and 
conjugacy problems in the braid groups},   Advances  Math.  139(1998), 322-353.

\bibitem{Bri0}
M.G.Brin, {\em  The chameleon groups of Richard J. Thompson: 
automorphisms and dynamics}, 
Inst. Hautes \'Etudes Sci. Publ. Math. No. 84, 1996, 5-33.

\bibitem{Bri}
M.G.Brin, {\em The Algebra of Strand Splitting. I. A Braided Version of Thompson's Group V}, math.GR/0406042.  
\bibitem{Bri2}
M.G.Brin, {\em The Algebra of Strand Splitting.II. A Presentation for the Braid Group on One Strand}, math.GR/0406043. 




\bibitem{br}       
K.S.Brown, {\em  Presentations for groups acting on       
simply-connected complexes}, J. Pure Appl. Algebra 32(1984),  1-10.      
     
\bibitem{BR} K.S.Brown, {\em Finiteness properties of groups}, Proceedings of the     
Northwestern conference on cohomology of groups (Evanston, Ill., 1985), J. Pure Appl. Algebra 44     
(1987), no. 1-3, 45--75.       
      
\bibitem{br2}       
K. S. Brown, {\em The Geometry of Finitely Presented Infinite Simple       
Groups}, Algorithms and Classification in Combinatorial Group Theory       
(G. Baumslag and C. F. Miller III , eds), MSRI Publications,       
vol. 23. Springer-Verlag (Berlin, Heidelberg, New-York), 1992, 121-136.        

\bibitem{Ca}
D. Calegari, {\em Circular groups, Planar groups and the Euler class}, 
Proceedings of the Casson Fest, Geom. Topol. Monogr. 7(2004), 431-491. 

      
\bibitem{CFP}      
J.W.Cannon, W.J. Floyd, and W.R. Parry,       
{\em Introductory notes on Richard Thompson's groups},       
Enseign. Math. 42(1996), 215-256.  
 
      
\bibitem{De1}
P.Dehornoy, {\em Geometric presentations for Thompson's groups}, 
math.GR/0407096. 

\bibitem{De2}
P.Dehornoy, {\em The group of parenthesized braids}, 
math.GR/0407097. 
      

\bibitem{DDRW}
P.Dehornoy, I.Dynnikov, Ivan, D.Rolfsen and B.Wiest, 
{ Why are braids orderable?}, 
{\em  Panoramas et Synth\`eses},  SMF, Paris, 2002.
 



\bibitem{Dy}
I.A.Dynnikov,
{\em Three-page representation of links}, 
Uspekhi Mat. Nauk 53 (1998), 237-238; 
translation in Russian Math. Surveys 53(1998),  1091-1092.    
     
   

\bibitem{ECHLPT} 
D.B.A.Epstein, J.W.Cannon, D.F. Holt, S.V.F. Levy, M.S.Paterson
and  W.P.Thurston, {Word processing in groups}, {\em 
Jones and Bartlett Publishers}, Boston, MA, 1992.  
  

\bibitem{Fa}
D.S.Farley, {\em  Finiteness and $\rm CAT(0)$ 
properties of diagram groups},   Topology  42(2003),  1065-1082. 


   
\bibitem{FG}      
L.Funar and R.Gelca, {\em On the groupoid of transformations of rigid      
  structures on surfaces},       
J.Math.Sci.Univ.Tokyo 6(1999), 599-646.      
   
\bibitem{FK} 
L.Funar and C.Kapoudjian, {\em On a universal mapping class group  
in genus zero}, G.A.F.A. 14(2004), 965-1012. 



\bibitem{GeSh}
S.M.Gersten and H. Short, {\em 
Small cancellation theory and automatic groups II},
Invent. Math. 105(1991), 641-662.

    
\bibitem{GS}      
E.Ghys and V.Sergiescu, {\em Sur un groupe remarquable de      
  diff\'eomorphismes du cercle}, Comment.Math.Helvetici 62(1987), 185-239.      
 

\bibitem{GrS}
P.Greenberg and V.Sergiescu, {\em 
An acyclic extension of the braid group}, 
Comment. Math. Helv.  66(1991),  109-138. 


\bibitem{GN}
R.I.Grigorchuk, V.V.Nekrashevich, V.I.Sushchanskii, {\em 
Automata, dynamical systems and infinite groups,} 
Proc. Steklov Inst. Math. 231(2000), 134-214. 


\bibitem{Guba}
V.Guba, {\em Polynomial isoperimetric inequalities for 
Richard Thompson's groups $F$, $T$, and $V$}, 
Algorithmic problems in groups and semigroups (Lincoln, NE, 1998),  
91-120, {\em Trends Math., Birkh\"auser Boston,} Boston, MA, 2000.

\bibitem{GuS}
V. Guba and M. Sapir,{\em  Diagram groups}, 
Mem. Amer. Math. Soc.  130(1997),  no. 620.
     


\bibitem{HT}      
A. Hatcher and W. Thurston, {\em  A presentation for the mapping class      
  group of a closed orientable surface}, Topology 19(1980), 221-237.      
      
      
 
\bibitem{KS}       
C. Kapoudjian and V.Sergiescu,
{\em An extension of the Burau representation to a mapping class group 
associated to Thompson's group $T$}, math.GT/0302300. 
   
\bibitem{KoS}
M.Kontsevich, Y.Soibelman, {\em Affine structures and non-archimedean 
spaces}, math.AG/0406564.

   
    

     
\bibitem{LS}      
P. Lochak and L. Schneps,      
\newblock {\em On universal Ptolemy-Teichmuller groupoid},       
\newblock { in  Geometric Galois Theory}, L.M.S. Lecture      
Notes Ser.,243, Cambridge Univ.Press, 1997.       
      
 
 

\bibitem{Mas}
D.R.Mason, {\em  On the $2$-generation of 
certain finitely presented infinite simple groups}, 
J. London Math. Soc. (2) 16(1977), 229--231.
   

\bibitem{Mo}
L.Mosher, {\em  Mapping class groups are automatic}, 
Ann. of Math. (2) 142(1995), 303-384.


  
     
      
\bibitem{pe}       
R.C.Penner, {\em The universal Ptolemy group and its completions},   
Geometric Galois actions, 2, 293--312, L.M.S. Lecture Notes Ser., 243,
Cambridge Univ. Press, Cambridge, 1997.      
 
     
    

 
\bibitem{Se} 
V.Sergiescu, {\em Graphes planaires et pr\'esentations des groupes de 
tresses}, Math. Zeitsch. 214(1993), 477-490. 
  


\bibitem{Wi}
B.Wiest, {\em Diagram groups, braid groups, and orderability}, 
J. Knot Theory Ramifications  12(2003),   321--332. 
 
\bibitem{Wit}
D. Witte, {\em Arithmetic groups of higher $Q$-rank cannot act on $1$-manifolds},   Proc. Amer. Math. Soc.  122(1994),   333--340.
      
\end{thebibliography}

}

\end{document}